\documentclass[a4paper,11pt]{amsart}
\usepackage[T1]{fontenc}
\usepackage{amsmath,amssymb,mathtools,microtype,mathrsfs,tikz}
\usepackage{lmodern}
\usepackage{amsmath}

\newcommand{\colhead}[3]{%
 \mathord{\ooalign{%
  \hfil$#1$\hfil\cr
  \hfil$\hphantom{#2}$\hfil\cr
  \hfil$\hphantom{#3}$\hfil\cr}}}

\newcommand{\blockmat}[4]{%
 \overset{
  \begin{matrix}
   \colhead{r}{#1}{#3}&\colhead{q}{#2}{#4}
  \end{matrix}
 }{
  \begin{bmatrix}
   #1&#2\\
   #3&#4
  \end{bmatrix}
 }
 \;\begin{matrix}r\\q\end{matrix}}
\usepackage[hidelinks]{hyperref}
\hypersetup{pdftitle={A subharmonicity principle with applications to product inequalities}}
\theoremstyle{plain}
\newtheorem{theorem}{Theorem}
\newtheorem{proposition}{Proposition}
\newtheorem*{conjecture}{WPI Conjecture}
\newtheorem{corollary}{Corollary}
\newtheorem{lemma}{Lemma}[section]
\theoremstyle{definition}
\newtheorem{example}{Example}[section]
\newtheorem{question}{Question}[section]
\numberwithin{equation}{section}
\newcommand{\E}{\mathbb E}
\newcommand{\R}{\mathbb R}
\newcommand{\C}{\mathbb C}
\newcommand{\HH}{\mathbb H}
\newcommand{\F}{\mathbb F}
\newcommand{\PP}{\mathcal P}
\newcommand{\W}{\mathcal W}
\newcommand{\dd}{\mathrm d}
\newcommand{\tr}{\operatorname{Re}\operatorname{tr}}

\DeclareMathOperator{\diag}{diag}
\DeclareMathOperator{\rank}{rank}
\DeclareMathOperator{\Cov}{Cov}
\title[Subharmonicity principle and product inequalities]{A subharmonicity principle with applications to product inequalities}
\author[Haoming Wang]{Haoming Wang}
\address{Center for Combinatorics, Nankai University, Tientsin {\rm300071}, China}
\email{wanghm37@nankai.edu.cn}
\date{}
\subjclass[2020]{Primary 60E15; Secondary 60G44, 60J60, 62H10}
\keywords{Wishart product inequality, logarithmic subharmonicity, Brownian motion, squared Bessel process, submartingale}
\begin{document}
\begin{abstract}
This paper establishes a subharmonicity principle with a probabilistic formulation in terms of commuting Markov operators and hitting times. The resulting monotonicity of conditional product moments yields a proof of the real Wishart product inequality conjecture in arbitrary dimension, as well as product inequalities for complex and quaternionic Gaussian and Wishart distributions. Further functional product inequalities and non-Gaussian extensions are also obtained.
\end{abstract}
\maketitle

\section{Introduction}
The purpose of this paper is to identify a potential-theoretic mechanism underlying product inequalities for Gaussian and Wishart distributions over the reals, complexes and quaternions. Our starting point is a subharmonicity principle
formulated in terms of commuting Markov transition operators and established by applying Doob's optional sampling theorem to submartingales at hitting times of prescribed levels. The applicability of this principle to product inequalities is governed by the compatibility of logarithmic and fundamental potentials, for which dimension two is critical. 

In \cite{GenestOuimetRichards2024}, Genest, Ouimet and Richards
formulated the following conjecture concerning joint moments
of disjoint principal minors of Wishart matrices.

\begin{conjecture}[{\cite[Conjecture~1.1]{GenestOuimetRichards2024}}]
\label{conj:wpi}
Fix a positive integer $p$. Let $S\sim \mathcal W_p^{(1)}(n,\Sigma)$, where
$\Sigma>0$ and $n>(p-1)/2$.
Partition both $S$ and $\Sigma$ into blocks of
sizes $p_1,\ldots,p_m$, where $p_1+\cdots+p_m=p$.
For any $\nu_1,\ldots,\nu_m\geq0$,
\begin{equation}\label{eq:wpi}
    \E\prod_{j=1}^m\det(S_{jj})^{\nu_j}
    \geq
    \prod_{j=1}^m\E\det(S_{jj})^{\nu_j}.
\end{equation}
\end{conjecture}
Here and throughout, the superscripts $(1)$, $(2)$ and
$(4)$ indicate the real, complex and quaternionic cases, respectively. The notation $A>0$ (resp. $A\geq0$) means that $A$
is positive definite (resp. positive semidefinite). The law $S\sim\W_p^{(\beta)}(n,\Sigma), \beta \in \{1,2,4\}$ is characterized by its Laplace transform 
\begin{equation}
 \E e^{-\tr(TS)} = \det\left(I_p+\frac{2}{\beta} \Sigma^{1/2} T\Sigma^{1/2}\right)^{-n}\ (T\geq0).
\end{equation}
The determinant of a Hermitian matrix is understood as the product of its eigenvalues, also in the noncommutative case.

In a subsequent paper, Genest, Ouimet and Richards
\cite[Theorem~3.1]{GenestOuimetRichards2025} derived an explicit formula for the joint moment of two
disjoint principal minors in terms of the Gaussian
hypergeometric function of matrix argument.
Their Corollary~3.4 uses this formula to give a quantitative version
of the WPI Conjecture for $m=2$.
Our first result settles this conjecture in arbitrary
dimension and characterizes the equality cases. 

\begin{theorem}\label{thm:real}
The WPI Conjecture holds.
Moreover, inequality~\eqref{eq:wpi} remains valid for $\Sigma\ge 0$ and $n\ge (p-1)/2$.
If $\Sigma_{jj}>0$ whenever $\nu_j>0$, equality holds
if and only if $\Sigma_{ij}=0$ for all $i\ne j$
with $\nu_i\nu_j\ne 0$.
\end{theorem}

In addition to settling the conjecture,
Theorem~\ref{thm:real} extends WPI to the endpoint of the continuous shape range and to positive semidefinite
covariances. In the real case, the remaining admissible shapes require a different argument, which is given in Section~\ref{sec:singular-real}.
In the complex and quaternionic cases, however, the subharmonicity method covers the entire admissible shape range.

For $\beta \in \{1,2,4\}$ and any covariance parameter $\Sigma>0$,
the admissible shape parameters of the Wishart distribution
$\mathcal W_p^{(\beta)}(n,\Sigma)$ constitute the Gindikin
set $\beta\mathscr A_p$ (see \cite[Theorem~VII.3.1]{FarautKoranyi1994}), where 
\begin{equation}\label{eq:gindikin}
     \mathscr A_p=
 \left\{\frac{k}{2}\mathrel{:} k=0,\ldots,p-2\right\}
 \cup
 \left[\frac{p-1}{2},\infty\right).
\end{equation}

\begin{corollary}\label{cor:ch}
Let $S\sim \mathcal W_p^{(\beta)}(n,\Sigma)$, where
$\beta=2$ or $4$, $\Sigma\geq0$ and
$n\in\beta \mathscr A_{p}$.
Then inequality~\eqref{eq:wpi} holds.
If $\Sigma_{jj}>0$ and $2n>\beta(p_j-1)$
whenever $\nu_j>0$, equality holds if and only if
$\Sigma_{ij}=0$ for all $i\ne j$ with $\nu_i\nu_j\ne0$.
\end{corollary}

The different shape ranges covered by our argument
can be understood in terms of logarithmic potentials.
A product of nonnegative submartingales need not be
a submartingale, whereas the convexity and monotonicity of the exponential function imply that if the logarithms of positive processes are submartingales, then any integrable product of positive powers of those processes is a submartingale.  

To explain the role of dimension, let
\begin{equation}\label{eq:fundamental-potentials}
 \Phi_d(r)=
 \begin{cases}
 r^{2-d}/(2-d),&d\ne2,\\
 \log r,&d=2,
 \end{cases}
\end{equation}
be the increasing radial fundamental potential of the
$d$-dimensional Laplacian. The function $\log r$ is an
increasing convex function of $\Phi_d(r)$ if and only if
$d\geq2$. Thus dimension two is precisely the threshold
at which the fundamental potential is compatible with
the logarithmic potential.

For partitioned Wishart matrices, conditional regression leads to comparisons in the noncentrality parameter. Under a suitable Wishart transition in a direction
$H=cc^*$, the Schur complement of a fixed positive definite principal block evolves as a squared Bessel process. In the successive comparisons, the corresponding dimensions of squared Bessel processes are
$2n-\beta(p_1+p_j-2)$ (see Section~\ref{sec:bessel}).
The continuous shape range ensures that these
dimensions are at least two whenever $m\geq3$, as required for logarithmic subharmonicity.

At the discrete shape parameters $n=k\beta/2$,
with $k=0,\ldots,p-1$, a Wishart matrix admits the
representation $S=X^*X$, where $X$ has $k$ independent
centered Gaussian rows with common covariance parameter $\Sigma$.
If $\Sigma>0$, then $S$ has rank $k$ almost surely.
Varying a single row along a fixed direction reduces
the logarithm of a regularized determinant, up to an
additive constant, to
\begin{equation}\label{eq:regularized-potential}
 \psi_d(u) = \log(b+c|u-u_0|^2),\quad u\in  \mathbb R^d,
\end{equation}
where $b>0$, $c\geq0$ and $u_0\in\mathbb R^d$. For $c>0$, $\psi_d$  is subharmonic
if and only if $d\geq2$. Taking $d=\beta$, we see that the complex and quaternionic cases satisfy the dimension requirement for the comparison of a single Gaussian row.

The following quaternionic Gaussian product inequality
(GPI) follows from Corollary~\ref{cor:ch} by taking
$n=2$ and $p_1=\cdots=p_m=1$.

\begin{corollary}\label{thm:gauss}
Let $Z=(Z_1,\ldots,Z_p)$ be a centered quaternionic
Gaussian row of covariance $\Sigma\geq0$.
For any $\nu_1,\ldots,\nu_p>0$,
\begin{equation}\label{eq:gpi}
 \E\prod_{j=1}^p|Z_j|^{\nu_j}
 \geq
 \prod_{j=1}^p\E|Z_j|^{\nu_j}.
\end{equation}
If $\Sigma_{jj}>0$ for every $j$, equality holds
if and only if the coordinates of $Z$ are independent.
\end{corollary}

The real and complex counterparts of~\eqref{eq:gpi}
were established by Ouimet and Greaves~\cite{OuimetGreaves2026}
and Barthe and Cordero-Erausquin~\cite{BartheCordero2025},
respectively. Corollary~\ref{cor:ch} provides an alternative proof
in the complex case, together with a characterization
of equality. However, the real case is not covered by our submartingale method, since logarithmic subharmonicity fails in dimension one. This motivates several open questions in the Discussion.

The key step in our proofs is to establish the monotonicity
of conditional product moments in the Loewner order on
positive semidefinite cones, which, together with spectral association obtained from multivariate total positivity of order two (MTP$_2$),
yields~\eqref{eq:wpi}. Conditional regression reduces the required monotonicity to
comparisons in the noncentrality parameter.
The relevant dimension conditions then arise from the
dynamics of Schur complements under Wishart diffusions
and ensure the submartingale property of regularized
logarithmic determinants, thereby allowing us to apply the
subharmonicity principle. The same subharmonicity argument applies to shifted Gaussian rows over the complex numbers and the quaternions.

This paper is organized as follows.
Section~\ref{sec:association} reduces the product
inequality to the monotonicity of conditional moments,
which is further reduced to row comparisons through
conditional regression in Section~\ref{sec:conditional}.
Section~\ref{sec:potentials} develops the submartingale
monotonicity principle, studies logarithmic potentials
and establishes the comparison of noncentral Wishart moments.
These results together yield proofs of WPI and its corollaries in Section~\ref{sec:wishart}. Section~\ref{sec:non-gaussian} treats the remaining discrete shapes in the real case and presents further product inequalities.
The appendices collect the necessary background
on Markov processes and Wishart diffusions.

\section{\texorpdfstring{MTP$_2$}{MTP2} and successive separation}\label{sec:association}

Let $\mathbb F$ be one of $\mathbb R$, $\mathbb C$ or $\mathbb H$, and let $\beta=\dim_{\R}\F$. Use the block partition from the introduction. Since blocks with zero exponents can be omitted in \eqref{eq:wpi}, assume all the remaining exponents are positive. The case $m=1$ is immediate, so we take $m\geq2$.

We separate the principal determinants one block at a time.

\begin{lemma}\label{lem:marginal-separation}
	To establish \eqref{eq:wpi}, it suffices to show, for the given partition and each successive principal marginal, that
\begin{equation}\label{eq:sep}
 \E\prod_{j=1}^m\det(S_{jj})^{\nu_j}
 \geq \E\det(S_{11})^{\nu_1}\,
       \E\prod_{j=2}^m\det(S_{jj})^{\nu_j}.
\end{equation}
\end{lemma}
\begin{proof}[Proof of Lemma~\ref{lem:marginal-separation}]
The Laplace transform shows that each principal submatrix of $S$
has a Wishart distribution with the same shape parameter and the
corresponding principal submatrix of the covariance parameter.
Apply~\eqref{eq:sep} first to blocks $1,\ldots,m$, then to blocks
$2,\ldots,m$, and continue. At each step the marginal moments
already extracted are nonnegative. The last factor is a single
marginal moment, giving~\eqref{eq:wpi}.
\end{proof}

We therefore fix one partition and examine the difference between the
two sides of~\eqref{eq:sep}.

\subsection{The conditional expectation}

Put $r=p_1$, $A= S_{11}$ and $D=\Sigma_{11}$. Define
\begin{equation}\label{eq:conditional-moment}
 G(a)=\E\!\left[\prod_{j=2}^m \det(S_{jj})^{\nu_j}
                \,\middle|\,A=a\right].
\end{equation}
The difference between the two sides of~\eqref{eq:sep} is
$\Cov(\det(A)^{\nu_1},W)$ where $W=\prod_{j=2}^m \det(S_{jj})^{\nu_j}$. Since $\det(A)^{\nu_1}$ is a measurable function of $A$, the total expectation formula shows that
this difference equals $\Cov(\det(A)^{\nu_1},G(A))$.

Suppose that $D>0$ and $2n>\beta(r-1)$, so that $A$ has a density
on the positive definite cone. Write $\lambda=(\lambda_1,\ldots,\lambda_r)$
for the eigenvalues of $D^{-1/2}AD^{-1/2}$ in decreasing order $\lambda_1>\dots>\lambda_r$.
By unitary invariance of this matrix,
\[
 A\overset{\mathrm d}{=}D^{1/2}Q\diag(\lambda)Q^*D^{1/2},
\]
where $Q$ is Haar distributed on the unitary group over $\F$ and
independent of $\lambda$. Suppose for the moment that $G$ admits a continuous
version on the positive definite cone.
For $x\in(0,\infty)^r$, define
\begin{equation}\label{eq:haar-average}
	 h(x)=\E_Q G\bigl(D^{1/2}Q\diag(x)Q^*D^{1/2}\bigr).
\end{equation}
Then $\E[G(A)\mid\lambda]=h(\lambda)$ by the independence of $Q$ and $\lambda$. Put $w_1(x)=\left(\prod_{i=1}^r x_i\right)^{\nu_1}$. Since $\det(A)^{\nu_1}=(\det D)^{\nu_1}w_1(\lambda)$ is measurable with respect to $\lambda$, a second
conditioning gives
\begin{equation}\label{eq:haar-covariance}
 \Cov\bigl(\det(A)^{\nu_1},W\bigr)
 =(\det D)^{\nu_1}\Cov\bigl(w_1(\lambda),h(\lambda)\bigr).
\end{equation}
Thus the separation problem reduces to a covariance inequality for two functions of the ordered eigenvalues. We next establish the positive dependence property
needed for this comparison.

\subsection{\texorpdfstring{MTP$_2$}{MTP2} and association}

For $x,y\in\R^r$, write $x\wedge y$ and $x\vee y$ for their
coordinatewise minimum and maximum. A nonnegative function $f$ is
multivariate totally positive of order two, or $\mathrm{MTP}_2$, if
\begin{equation}\label{eq:mtp}
 f(x\wedge y)f(x\vee y)\geq f(x)f(y)
\end{equation}
for all $x,y$. A random vector $X$ is associated if
$\Cov(f(X),g(X))\geq0$ for all coordinatewise nondecreasing
functions $f$ and $g$ for which the covariance exists
\cite[Definition~1.1]{EsaryProschanWalkup1967}.
A random vector with an $\mathrm{MTP}_2$ density is associated
\cite{KarlinRinott1980}. We shall also use two elementary closure
properties. Independent associated vectors are jointly associated,
by conditioning on one vector and using association for the other.
A vector of coordinatewise nondecreasing functions of an associated vector is also associated, by composition.

Association gives a nonnegative covariance for
coordinatewise nondecreasing functions.
To characterize equality in the separation inequality,
we need a strict version of this property.

\begin{lemma}\label{lem:strict-association}
	Let $X$ be an associated random vector and let
	$\Omega\subset(0,\infty)^r$ be an open set such that
	$\mathbb P(X\in\Omega)=1$ and
	$\mathbb P(X\in U)>0$ for every nonempty open
	subset $U$ of $\Omega$. Suppose that $w\in C^1((0,\infty)^r)$ satisfies
	$\partial_i w>0$ for $i=1,\ldots,r$ and that $h$ is continuous
	and coordinatewise nondecreasing on $(0,\infty)^r$.
	If the covariance exists, then $\Cov(w(X),h(X))>0$ if and only if
	$h$ is nonconstant on $\Omega$.
\end{lemma}

\begin{proof}[Proof of Lemma~\ref{lem:strict-association}]The strict positivity of the partial derivatives of
$w$ allows us to perturb $w$ locally without losing
coordinatewise monotonicity.

	Fix $\psi\in C_c^\infty(\Omega)$ and extend it by zero to
	$(0,\infty)^r$. On the support of $\psi$, each $\partial_i w$
	has a positive lower bound and each $\partial_i\psi$
	is bounded. Thus $w\pm\varepsilon\psi$ are both
	coordinatewise nondecreasing for sufficiently small
	$\varepsilon>0$. Applying association with $h$ to
	these two functions gives
	\[
	\varepsilon\bigl|\Cov(\psi(X),h(X))\bigr|
	\leq \Cov(w(X),h(X)).
	\]
	
	Hence, if the covariance vanishes, $h(X)$ must be uncorrelated
	with every smooth test function supported in $\Omega$.
	Equivalently,
	\[
	\E[\psi(X)\{h(X)-\E h(X)\}]=0 \ (\psi\in C_c^\infty(\Omega)).
	\]
	We now use the continuity of $h$ and the support
	assumption on $X$ to show that $h$ is constant.
	
	If $h(x)\ne\E h(X)$ at some $x\in\Omega$,
	there is an open neighborhood of $x$ in $\Omega$
	on which $h-\E h(X)$ has a fixed strict sign.
	Choosing a nonnegative test function supported
	in this neighborhood and positive at $x$
	contradicts the preceding identity, since every
	nonempty open subset of $\Omega$ has positive
	probability under the law of $X$.
	Thus $h$ is constant on $\Omega$.
	The converse is immediate.
\end{proof}

Dykstra and Hewett~\cite{DykstraHewett1978} studied the positive dependence of Wishart eigenvalues using total positivity. We recall the $\mathrm{MTP}_2$ property of the ordered eigenvalues and derive a strict covariance inequality.

\begin{lemma}\label{lem:spectrum}
	Let $A\sim\W_r^{(\beta)}(n,D)$, where $D>0$ and
	$2n>\beta(r-1)$. The ordered eigenvalues
	$\lambda_1>\cdots>\lambda_r>0$ of $D^{-1/2}AD^{-1/2}$
	have an $\mathrm{MTP}_2$ density and are associated.
	Put $w(x)=(\prod_{i=1}^r x_i)^\nu$ with $\nu>0$.
	Then, for every continuous, coordinatewise nondecreasing
	function $h$ on $(0,\infty)^r$,
	\begin{equation}\label{eq:spectralcov}
		\Cov(w(\lambda),h(\lambda))>0
	\end{equation}
	  whenever the covariance exists and $h$ is nonconstant on the open set $\Omega=\{x\in\R^r:x_1>\cdots>x_r>0\}$.
\end{lemma}

\begin{proof}[Proof of Lemma~\ref{lem:spectrum}]
	By the change of variables, the joint density
	$f(\lambda)$ of the ordered eigenvalues $\lambda=(\lambda_1,\ldots,\lambda_r)$
	is proportional on $\Omega$ to
	\begin{equation}\label{eq:spectral-density}
		\underbrace{\prod_{i=1}^r
			\lambda_i^{n-1-\beta(r-1)/2}}_{(1)}
		\underbrace{\prod_{i=1}^r e^{-\beta\lambda_i/2}}_{(2)}
		\underbrace{\prod_{i<j}(\lambda_i-\lambda_j)^\beta}_{(3)}
	\end{equation}
	and vanishes outside $\Omega$.
	It therefore suffices to verify~\eqref{eq:mtp}
	for $x,y\in\Omega$. On this domain, $x\wedge y$ and $x\vee y$ also belong to
	$\Omega$.
	
	For factors (1) and (2), equality holds
	in~\eqref{eq:mtp}, since
	$(x_i\wedge y_i,x_i\vee y_i)$ is a rearrangement
	of $(x_i,y_i)$ for every $i$.
	It remains to verify the inequality for factor (3).
	For each $i<j$, put $d_{ij}(x)=x_i-x_j$.
	We need to show that
	\[
	d_{ij}(x\wedge y)d_{ij}(x\vee y)
	\geq d_{ij}(x)d_{ij}(y).
	\]
	If $x_i\leq y_i$ and $x_j\leq y_j$, or if both
	inequalities are reversed, equality holds.
	Otherwise, exchanging $x$ and $y$ if necessary,
	we may assume that $x_i\leq y_i$ and $y_j\leq x_j$.
	In this case,
	\[
        \begin{aligned}
         &d_{ij}(x\wedge y)d_{ij}(x\vee y)
           -d_{ij}(x)d_{ij}(y)\\
         &=(x_i-y_j)(y_i-x_j)
           -(x_i-x_j)(y_i-y_j)\\
         &=(y_i-x_i)(x_j-y_j)\geq0.
        \end{aligned}
	\]
	Raising both sides of these inequalities to the power $\beta$
	and multiplying over $i<j$ gives~\eqref{eq:mtp}
	for factor (3).
	Thus $f$ is $\mathrm{MTP}_2$, and the ordered
	eigenvalues are associated.
	
	Since $f$ is strictly positive on $\Omega$,
	every nonempty open subset of $\Omega$ has positive
	probability. Moreover, $\partial_i w(x)=\nu w(x)/x_i>0$
	for every $i$. Lemma~\ref{lem:strict-association}
	therefore gives~\eqref{eq:spectralcov}.
\end{proof}

\subsection{A sufficient condition for separation}

For Hermitian matrices, $a_0\leq a_1$ in the Loewner order means
that $a_1-a_0$ is positive semidefinite. Combining spectral association with Haar averaging gives a sufficient condition for the separation inequality.

\begin{lemma}\label{lem:spectral-reduction}
	Let $A\sim\W_r^{(\beta)}(n,D)$, where $D>0$ and
	$2n>\beta(r-1)$, and let $W$ be a nonnegative and integrable random variable.
	Suppose that $G(a)=\E[W\mid A=a]$ admits a continuous,
	Loewner-nondecreasing version on the positive definite cone.
	Then, for $\nu_1>0$,
	\[
	\Cov(\det(A)^{\nu_1},W)\geq0
	\]
	whenever the covariance exists, with equality if and
	only if this version of $G$ is constant on the positive definite cone.
\end{lemma}

\begin{proof}[Proof of Lemma~\ref{lem:spectral-reduction}]
	Let $h$ be the average of $G$ defined
	in~\eqref{eq:haar-average}.
	The integrand in~\eqref{eq:haar-average} is jointly
	continuous in $(x,Q)$ and hence locally bounded in
	$x$, uniformly in $Q$, by compactness of the unitary
	group. Dominated convergence implies that $h$ is continuous.
	For $x,y\in(0,\infty)^r$ with $x_i\leq y_i$ for
	every $i$, we have
	\[
	D^{1/2}Q\diag(x)Q^*D^{1/2}
	\leq D^{1/2}Q\diag(y)Q^*D^{1/2}
	\]
	in the Loewner order for every unitary matrix $Q$.
	By the assumed monotonicity of $G$, taking expectations
	over $Q$ gives $h(x)\leq h(y)$.
	Thus $h$ is coordinatewise nondecreasing.
	
	By~\eqref{eq:haar-covariance} and
	Lemma~\ref{lem:spectrum}, we obtain
	\[
	\Cov(\det(A)^{\nu_1},W)
	=(\det D)^{\nu_1}\Cov(w_1(\lambda),h(\lambda))
	\geq0.
	\]
	
	Suppose that equality holds. By Lemma~\ref{lem:spectrum},
	$h$ is constant on $\Omega$. Since $(t,\ldots,t)$
	belongs to the closure of $\Omega$ for every $t>0$,
	continuity implies that $h(t,\ldots,t)$ takes the same value for all $t>0$.
	By~\eqref{eq:haar-average},
	\[
	h(t,\ldots,t)=\E_Q G(tD^{1/2}QQ^*D^{1/2})
	=G(tD).
	\]
	Hence $G$ is constant along the ray $\{tD:t>0\}$.
	
	For any positive definite matrix $a$, there exist
	$0<s\leq t$ such that $sD\leq a\leq tD$.
	The Loewner monotonicity of $G$ gives
	$G(sD)\leq G(a)\leq G(tD)$.
	Since the two bounds are equal, $G$ is constant
	on the positive definite cone.
	Conversely, if $G$ is constant, then so is $h$,
	and equality follows from~\eqref{eq:haar-covariance}.
\end{proof}

Under the same assumptions on $G$, the covariance
inequality in Lemma~\ref{lem:spectral-reduction}
also holds for any almost surely positive definite
random matrix $A$ such that $D^{-1/2}AD^{-1/2}$,
for some $D>0$, is unitarily invariant and has associated ordered eigenvalues.
In this inequality, the determinant power may also be replaced by any
nondecreasing function of the determinant, provided
that the covariance exists. The remaining task is to construct a continuous, Loewner-nondecreasing version of $G$ and determine when $t\mapsto G(tD)$ is nonconstant.

\section{Conditional regression and row comparisons}
\label{sec:conditional}

Let $\PP_d$ be the cone of positive semidefinite
Hermitian $d\times d$ matrices over $\F$, where
$\beta=\dim_{\R}\F\in\{1,2,4\}$.
Throughout this section, suppose that $\Sigma\geq0$,
$\Sigma_{11}>0$ and $\nu_j>0$ for $j=2,\ldots,m$.
Put $r=p_1$ and $q=p-r$, and partition $S$ and $\Sigma$ as
\begin{equation}\label{eq:confomal-partition}
 S=\blockmat{A}{S_{AB}}{S_{BA}}{B},
 \qquad
 \Sigma=\blockmat{\Sigma_{AA}}{\Sigma_{AB}}
                 {\Sigma_{BA}}{\Sigma_{BB}},
\end{equation}
where $A=S_{11}$ and $D=\Sigma_{AA}=\Sigma_{11}$. Define $C=\Sigma_{AA}^{-1}\Sigma_{AB}$ and $\Psi=\Sigma_{BB}-\Sigma_{BA}C$.

We reduce the problem of proving that $G$, defined in~\eqref{eq:conditional-moment}, is nondecreasing in the Loewner order to a sequence of comparisons involving individual rows.

\subsection{Conditional Wishart regression} We first recall a conditional representation of partitioned Wishart matrices. We use the convention that a centered $\F$-Gaussian row of covariance $\Psi$ has the distribution of $\beta^{-1/2}\xi\Psi^{1/2}$, where the real components of $\xi$ are independent standard normal variables.

\begin{lemma}[Conditional Wishart regression]\label{lem:cond-regression}
Suppose that $2n>\beta(p-1)$ or that $n=k\beta/2$
for an integer $k\geq r$.
For $a>0$, let $L\in\F^{r\times r}$ satisfy $L^*L=a$.
A version of the conditional distribution of $B$
given $A=a$ is
\begin{equation}\label{eq:cond-regression}
 B\mid A=a\ \overset{\mathrm d}{=}
 Y+(Z+LC)^*(Z+LC),
\end{equation}
where $Z$ has $r$ independent centered Gaussian rows
of covariance $\Psi$, and
$Y\sim\W_q^{(\beta)}(n-\beta r/2,\Psi)$
is independent of $Z$.
\end{lemma}

\begin{proof}[Proof of Lemma~\ref{lem:cond-regression}]
Let $\widetilde S$ be a Wishart matrix with shape $n$ and covariance parameter $\diag(D,I_q)$, and partition it as in~\eqref{eq:confomal-partition}.
We first show that
\begin{equation}\label{eq:standard-regression}
 \widetilde S\mid\widetilde A=a\ \overset{\mathrm d}{=}\
 \begin{bmatrix}a&L^*Z_0\\ Z_0^*L&Y_0+Z_0^*Z_0\end{bmatrix},
\end{equation}
where $Z_0$ has $r$ independent centered Gaussian rows of covariance
$I_q$ and is independent of $Y_0\sim\W_q^{(\beta)}(n-\beta r/2,I_q)$.

For $2n>\beta(p-1)$, this is the standard conditional Wishart
decomposition~\cite[Theorem~7(i)--(ii)]{MassamNeher1997}. For discrete shape parameters, the same decomposition follows from the Gaussian
representation of Wishart matrices.

Let $n=k\beta/2$, and write $\widetilde S=(X_A,E)^*(X_A,E)$, where
$X_A$ and $E$ are independent centered Gaussian matrices
with $k$ rows and respective row covariances
$\Sigma_{AA}$ and $I_q$.
Conditionally on $X_A^*X_A=a$, choose a unitary matrix $Q=Q(X_A)$
such that
\[
 QX_A=\begin{bmatrix}L\\0\end{bmatrix},
 \quad
 QE=\begin{bmatrix}Z_0\\E_0\end{bmatrix}.
\]
By left unitary invariance, $QE$ has the same
distribution as $E$ and is independent of $X_A$.
Taking $Y_0=E_0^*E_0$ gives
\eqref{eq:standard-regression}.

For the stated covariance $\Sigma$, put
\[
 T=\begin{bmatrix}I_r&C\\0&\Psi^{1/2}\end{bmatrix}.
\]
Then $T^*\diag(D,I_q)T=\Sigma$, so the Wishart transform gives $S\overset{\mathrm d}{=}T^*\widetilde S T$ even when $T$ is singular. This congruence leaves the first block unchanged. Since the matrix on the right of
\eqref{eq:standard-regression} equals
$\diag(0,Y_0)+(L,Z_0)^*(L,Z_0)$, its transformed complementary block is
\[
 \Psi^{1/2}Y_0\Psi^{1/2}
 +(Z_0\Psi^{1/2}+LC)^*(Z_0\Psi^{1/2}+LC).
\]
Thus $Y=\Psi^{1/2}Y_0\Psi^{1/2}$ and $Z=Z_0\Psi^{1/2}$ have the
required laws and are independent.
\end{proof}


For $U\in\PP_q$, label its diagonal blocks by $2,\ldots,m$
in the order inherited from $B$, and put
\begin{equation}\label{eq:remaining-determinant-product}
F(U)=\prod_{j=2}^m\det(U_{jj})^{\nu_j},
\end{equation}
where $U_{22},\ldots,U_{mm}$ denote the principal blocks of $U$
corresponding to those of $B$.
By Lemma~\ref{lem:cond-regression}, the conditional
moment can be written as
\begin{equation}\label{eq:conditional-reduction}
 G(a)=\E F(Y+(Z+LC)^*(Z+LC)),
 \quad L^*L=a.
\end{equation}
With $L=a^{1/2}$ in~\eqref{eq:conditional-reduction},
the integrand is continuous in $a$.
For $a$ in a compact subset of the positive
definite cone, the integrand is bounded uniformly in $a$ by a fixed polynomial in $1+\tr Y+\tr(Z^*Z)$. This bound is integrable since $Y$ and $Z$ have moments of all orders. Thus $G$ is continuous by dominated convergence.

\subsection{Reduction to row comparisons}

Although we have used $L=a^{1/2}$ to establish
continuity, the representation~\eqref{eq:conditional-reduction}
allows any square factor $L$ satisfying $L^*L=a$.
We shall use this freedom to compare two conditional
moments by varying one row at a time.

If $M\in\PP_q$ is the contribution of the rows held fixed, the
corresponding function of the remaining matrix is
\begin{equation}\label{eq:cond-observable}
 F_M(U)=F(M+U)=\prod_{j=2}^m\det(M_{jj}+U_{jj})^{\nu_j}.
\end{equation}

\begin{lemma}
\label{lem:row-reduction}
Suppose that the assumptions of
Lemma~\ref{lem:cond-regression} hold.
Let $g$ be a centered Gaussian row of covariance
$\Psi$, independent of $Y$.
If, for every $M\in\PP_q$ and every deterministic
row $v\in\F^{1\times q}$, the function
\begin{equation}\label{eq:row-comparison}
 t\mapsto\E F(M+Y+(g+tv)^*(g+tv))
\end{equation}
is nondecreasing on $[0,\infty)$, then $G$ is
nondecreasing with respect to the Loewner order
on the positive definite cone.
\end{lemma}

\begin{proof}[Proof of Lemma~\ref{lem:row-reduction}]
Let $0<a_0\leq a_1$. Choose a unitary matrix $Q$
such that
\[
 a_0^{-1/2}a_1a_0^{-1/2}
 =Q\diag(t_1^2,\ldots,t_r^2)Q^*,
 \qquad t_i\geq1.
\]
Put $R=Q^*a_0^{1/2}$. Then $a_0=R^*R $ and $ a_1=R^*\diag(t_1^2,\ldots,t_r^2)R.$

This representation allows us to compare the conditional
moments by varying one row at a time.
For $i=0,\ldots,r$, define
\[
 L^{(i)}=\diag(t_1,\ldots,t_i,1,\ldots,1)R,
 \qquad a^{(i)}=(L^{(i)})^*L^{(i)}.
\]
Since $a^{(0)}=a_0$ and $a^{(r)}=a_1$, it suffices
to establish $G(a^{(i-1)})\leq G(a^{(i)})$
for $i=1,\ldots,r$.

Fix $i$ and write $z_j$ and $v_j$ for the $j$th
rows of $Z$ and $RC$, respectively.
In passing from $L^{(i-1)}$ to $L^{(i)}$,
the $i$th row of the shifted Gaussian matrix
changes from $z_i+v_i$ to $z_i+t_iv_i$.
The contribution of the other rows is
\[
 \begin{aligned}
 M={}&\sum_{j=1}^{i-1}(z_j+t_jv_j)^*(z_j+t_jv_j)+\sum_{j=i+1}^{r}(z_j+v_j)^*(z_j+v_j).
 \end{aligned}
\]
Putting $g=z_i$ and $v=v_i$, we obtain
from~\eqref{eq:conditional-reduction}
\[
 \begin{aligned}
 G(a^{(i-1)})&=\E F(M+Y+(g+v)^*(g+v)),\\
 G(a^{(i)})&=\E F(M+Y+(g+t_iv)^*(g+t_iv)).
 \end{aligned}
\]

Since $M$ is independent of $(Y,g)$, we may first
take expectations over $Y$ and $g$ with $M$ fixed.
By the assumed monotonicity in~\eqref{eq:row-comparison},
the inner expectation at $t=t_i$ is at least its value at $t=1$.
Taking expectations over $M$ gives
$G(a^{(i-1)})\leq G(a^{(i)})$.
Repeating the comparison for each row yields
$G(a_0)\leq G(a_1)$.
\end{proof}

Since $Y$ is independent of $g$, we have
\[
 \E F_M(Y+(g+tv)^*(g+tv))
 =\E_Y\bigl[
 \E_g F_{M+Y}((g+tv)^*(g+tv))
 \bigr].
\]
Consequently, to prove the monotonicity in~\eqref{eq:row-comparison}, it suffices to show that the function
\begin{equation}\label{eq:direct-row-comparison}
 t\mapsto\E F_M((g+tv)^*(g+tv))
\end{equation}
is nondecreasing for every $M\in\PP_q$ and every deterministic row $v$.


\begin{lemma}\label{lem:row-transform}
Suppose that $2n\geq\beta(p-1)$ and put $\alpha=n-\beta(r-1)/2$.
For $T>0$, write $J=(T^{-1}+2\Psi/\beta)^{-1}$. With $Y$ and $g$
as in Lemma~\ref{lem:row-reduction}, one has
\begin{equation}\label{eq:shifted-row-transform}
 \E e^{-\tr\{T[Y+(g+tv)^*(g+tv)]\}}
 =\left(\frac{\det J}{\det T}\right)^\alpha
 \exp\{-t^2vJv^*\}.
\end{equation}
\end{lemma}

\begin{proof}[Proof of Lemma~\ref{lem:row-transform}]
The Laplace transform of $Y$ is
$(\det J/\det T)^{n-\beta r/2}$.
For $\Psi>0$, completing the square in the Gaussian density gives
the Laplace transform of $(g+tv)^*(g+tv)$ as
$(\det J/\det T)^{\beta/2}\exp\{-t^2vJv^*\}$.
Independence of $Y$ and $g$ now gives \eqref{eq:shifted-row-transform}.
The formula extends to singular $\Psi$ by continuity. At the
endpoint $2n=\beta(p-1)$, Lemma~\ref{lem:cond-regression} applies
with $k=p-1\geq r$.
\end{proof}

In~\eqref{eq:shifted-row-transform}, $t$ enters only
through the noncentrality matrix $t^2v^*v$. It therefore suffices to compare the expectations of $F_M$ along this ray for each fixed $M$. 

Both this comparison and the direct row
comparison~\eqref{eq:direct-row-comparison}
can be established through subharmonicity,
using Wishart diffusions for the former and
Brownian motion for the latter.
We first develop the underlying monotonicity principle.

\section{Subharmonicity and submartingale monotonicity}
\label{sec:potentials}

To establish the comparisons  formulated  in the
preceding section, we seek a principle that converts
subharmonicity into monotonicity along prescribed rays.
The main ingredients are commutation of Markov
operators and optional sampling at hitting times.

In the following, we consider a right-continuous strong
Markov process $X_t \ (t\ge0)$ on $E$ with conservative transition function $P(t,x,\varGamma) \  (t\ge 0,x\in E,\varGamma\in\mathscr{B})$ (see Appendix~\ref{app:markov-process}).
For each initial state $x\in E$, the filtered probability
space under $\mathbb P_x$ is assumed to satisfy the
natural conditions.
\subsection{A monotonicity principle}
A nonnegative measurable function $f$ on $E$ is {subharmonic} with respect to the transition semigroup $P_t \ (t\ge0)$ if $P_tf\geq f$ for all $t\geq0.$ 
\begin{lemma}\label{lem:submartingale}
    For every $x\in E$ and every subharmonic function $f$ satisfying $\E_x f(X_t)<\infty$ for all $t\geq0$, the process $f(X_t) $ is a submartingale under $\mathbb P_x$.
\end{lemma}

We observe that the class of subharmonic functions
is preserved by Markov operators commuting with the
transition semigroup.

\begin{lemma}\label{lem:commutation}
Let $K(x,\varGamma)\ (x\in E,\varGamma\in\mathscr B) $ be a Markov kernel whose associated operator $K$ commutes with the transition semigroup $P_t \ ({t\geq0})$ so that
$KP_t=P_tK$ for every $t\geq0$.
If $f$ is a nonnegative subharmonic function, then $Kf$
is also a subharmonic function.
\end{lemma}

\begin{lemma}[Doob's optional sampling]\label{lem:optional}
Let $Y_t \ (t\ge0)$ be a right-continuous submartingale.
If $\sigma$ and $\tau$ are bounded stopping times
with $\sigma\leq\tau$, then $\E[Y_\tau\mid\mathscr F_\sigma]\geq Y_\sigma \textnormal{ a.s.}$
If $\tau$ is almost surely finite and the family
$\{Y_{t\wedge\tau}:t\geq0\}$ is uniformly integrable,
then $\E Y_\tau\geq\E Y_0$.
\end{lemma}
The proofs of Lemma \ref{lem:submartingale} and Lemma \ref{lem:commutation} are direct, and the proof of Lemma \ref{lem:optional} can be found in {\cite[Chapter~II,
Section~3]{RevuzYor1999}}. These three lemmas give the following principle.

\begin{proposition}[Monotonicity principle]
\label{thm:monotonicity}
Let $K$ satisfy the conditions of
Lemma~\ref{lem:commutation}, and let $f$ be a nonnegative subharmonic function.
Suppose that $Kf(X_t)$ is integrable for every $t \ge 0$ and
has right-continuous paths under $\mathbb P_x$.
If $\tau$ is an almost surely finite stopping time
and the family
\[
 \{Kf(X_{t\wedge\tau}):t\geq0\}
\]
is uniformly integrable under $\mathbb P_x$, then
\begin{equation}\label{eq:stopping-comparison}
 Kf(x)\leq\E_x Kf(X_\tau).
\end{equation}
In particular, if $X_\tau=y$ almost surely, then
$Kf(x)\leq Kf(y)$.
\end{proposition}
\begin{proof}[Proof of Proposition \ref{thm:monotonicity}]
By Lemmas~\ref{lem:submartingale}
and~\ref{lem:commutation}, the process
$Kf(X_t)$ is a submartingale under $\mathbb P_x$.
By assumption, it has right-continuous paths.  

By Lemma~\ref{lem:optional}, the almost sure finiteness
of $\tau$ and the uniform integrability of the
stopped family yield
\[
 Kf(x)=\E_x Kf(X_0)\leq\E_x Kf(X_\tau).
\]
If $X_\tau=y$ almost surely, this gives
$Kf(x)\leq Kf(y)$.
\end{proof}

For a process with continuous paths, a natural class
of stopping times is furnished by exit times. Let $D$
be an open subset of $E$ and define
\[
 \tau_D=\inf\{t\geq0:X_t\notin D\}.
\]
If $x\in D$ and $\tau_D$ is almost surely finite under
$\mathbb P_x$, the exit distribution
\[
 \omega_D(x,\varGamma)
 =\mathbb P_x(X_{\tau_D}\in\varGamma)
\]
is concentrated on $\partial D$. When $D$ is relatively compact and $g=Kf$ is
continuous on $\overline D$, the stopped family
is uniformly bounded. 
Under the
remaining hypotheses of Proposition~\ref{thm:monotonicity},
formula~\eqref{eq:stopping-comparison} becomes
\begin{equation}\label{eq:exit-comparison}
 g(x)\leq\int_{\partial D}g(y)\omega_D(x,\dd y).
\end{equation}
When the exit distribution is concentrated at a single
point, this inequality reduces to a pointwise comparison. We next consider such a situation on a cone.

Let $\mathcal C$ be a closed convex cone and let
$P_A \ {(A\in\mathcal C)}$ be a family of Markov
operators satisfying $P_0=I$ and $P_AP_B=P_{A+B}$ for all $A,B \in \mathcal C$.
We refer to Appendix~\ref{sec:wishart-semigroup}
for the corresponding definitions and examples of such operators.

For a nonzero $H\in\mathcal C$, the operators
$P_{tH}$ form a transition semigroup.
The ray $\{xH:x\geq0\}$ is invariant if a process
starting on this ray remains on it.

\begin{proposition}[Monotonicity along rays]\label{cor:ray}
Fix a nonzero $H\in\mathcal C$. Suppose that the ray $\{xH:x\geq0\}$
is invariant under $P_{tH} \ ({t\geq0})$ and that the restriction of
the semigroup to the ray admits a continuous Markov realization
$X_t^H=R_tH$. For every $0\leq x<y$, assume that
\[
 \mathbb P_{xH}(\tau_y<\infty)=1,\quad
 \tau_y=\inf\{t\geq0:R_t=y\}.
\]
Let $f$ be a nonnegative function that is subharmonic with respect to
$(P_{tH})_{t\geq0}$ and satisfies $P_Af(U)<\infty$ for all $A,U\in\mathcal C$.
For every $V\in\mathcal C$ such that $P_Vf$ is continuous on this ray,
\begin{equation}\label{eq:ray-monotonicity}
 P_Vf(xH)\leq P_Vf(yH)\ (0\leq x\leq y).
\end{equation}
\end{proposition}
    
\begin{proof}[Proof of Proposition \ref{cor:ray}]
By the semigroup property, $P_V$ commutes with $P_{tH}$
for every $t\geq0$. Moreover,
\[
 \E_{xH}P_Vf(R_tH)=P_{V+tH}f(xH)<\infty.
\]
The continuity of $P_Vf$ on the ray ensures that
$P_Vf(R_tH)\ (t\geq0)$ has continuous paths.

For $x<y$, the hitting time $\tau_y$ is the exit time
from $[0,y)H$ relative to the ray. Since
$R_{t\wedge\tau_y}\in[0,y]$, we have
\[
 \sup_{t\geq0}P_Vf(R_{t\wedge\tau_y}H)
 \leq\max_{0\leq z\leq y}P_Vf(zH)<\infty.
\]
The stopped family is therefore uniformly integrable.
As $R_{\tau_y}=y$ almost surely,
Proposition~\ref{thm:monotonicity} yields
\eqref{eq:ray-monotonicity}.
\end{proof}

\subsection{Logarithmic and fundamental potentials}

The applicability of Proposition~\ref{cor:ray}
requires both the finiteness of hitting times and the subharmonicity of the functions under consideration. As will be seen in the next subsection, both issues
are related to the properties of Bessel processes.

Let $R_t\ (t\geq0)$ be a Bessel process of dimension
$d>0$ and let $Z_t=R_t^2$ be the corresponding squared
Bessel process. On $(0,\infty)$ their respective
infinitesimal generators are
\[
 \mathcal A_d f(x)=\frac12 f''(x)+\frac{d-1}{2x}f'(x),
 \quad
 \mathcal B_d f(x)=2xf''(x)+df'(x),
\]
for test functions $f\in C^2(0,\infty)$.

\begin{lemma}\label{lem:besq-hitting}
Let $Z_t \ (t \ge 0)$ start at $x\geq0$. For $y>x$, let
$\tau_y=\inf\{t\geq0:Z_t=y\}$. Then $\E_x\tau_y={(y-x)}/{d}$.
In particular, $\tau_y$ is almost surely finite.
\end{lemma}

\begin{proof}[Proof of Lemma \ref{lem:besq-hitting}]
The squared Bessel equation
$\dd Z_t=2\sqrt{Z_t}\,\dd B_t+d\,\dd t$ gives
\[
 \E_xZ_{t\wedge\tau_y}
 =x+d\E_x(t\wedge\tau_y)\leq y.
\]
Since $t\wedge\tau_y\uparrow\tau_y$ as $t\to\infty$, monotone convergence yields
$\E_x\tau_y\leq(y-x)/d<\infty$.
Since $Z_{t\wedge\tau_y}\leq y$ and
$Z_{\tau_y}=y$ almost surely, bounded convergence
gives the asserted equality.
\end{proof}

We next consider the exit distribution of Bessel processes from
a bounded interval $(a,b)$ with $a>0$.
The fundamental potential $\Phi_d$ defined
in~\eqref{eq:fundamental-potentials} satisfies
$\Phi_d'(r)=r^{1-d}$ and $\mathcal A_d\Phi_d=0$.
It is therefore strictly increasing and harmonic
with respect to $\mathcal A_d$.

\begin{lemma}
\label{lem:scale}
For $0<a<r<b$, suppose that $R_0 = r$ and write
$T_c=\inf\{t\geq0:R_t=c\}$. Then
\begin{equation}\label{eq:annulus-exit}
 \mathbb P_r(T_b<T_a)
 =\frac{\Phi_d(r)-\Phi_d(a)}
        {\Phi_d(b)-\Phi_d(a)}.
\end{equation}
Consequently, for a continuous function $h$ on
$(0,\infty)$, the following conditions are equivalent.
\begin{enumerate}
\item[(i)] $h\circ\Phi_d^{-1}$ is convex on $\Phi_d((0,\infty))$.
\item[(ii)] $h(r)\leq\E_r h(R_{T_a\wedge T_b})$
for every $0<a<r<b$.
\end{enumerate}
If $h$ is twice continuously differentiable, they are also equivalent to
\begin{enumerate}
\item[(iii)] $\mathcal A_d h\geq0$ on $(0,\infty)$.
\end{enumerate}
\end{lemma}

\begin{proof}
Applying Lemma~\ref{lem:besq-hitting} to $R_t^2$
at level $b^2$, we see that $T_b$ is almost surely
finite. Since $\Phi_d$ is harmonic,
the stopped process
$\Phi_d(R_{t\wedge T_a\wedge T_b})$ is a bounded
martingale. The optional sampling lemma gives
\[
 \Phi_d(r)=\E_r\Phi_d(R_{T_a\wedge T_b}).
\]
As $R_{T_a\wedge T_b}$ takes values in $\{a,b\}$,
this identity yields~\eqref{eq:annulus-exit}.

The equivalence of (i) and (ii) follows from
\eqref{eq:annulus-exit} and the definition of convexity.
If $h=F\circ\Phi_d$ is twice continuously differentiable,
then
\[
 \mathcal A_d h
 =\frac12(\Phi_d')^2(F''\circ\Phi_d),
\]
which proves the last assertion.
\end{proof}

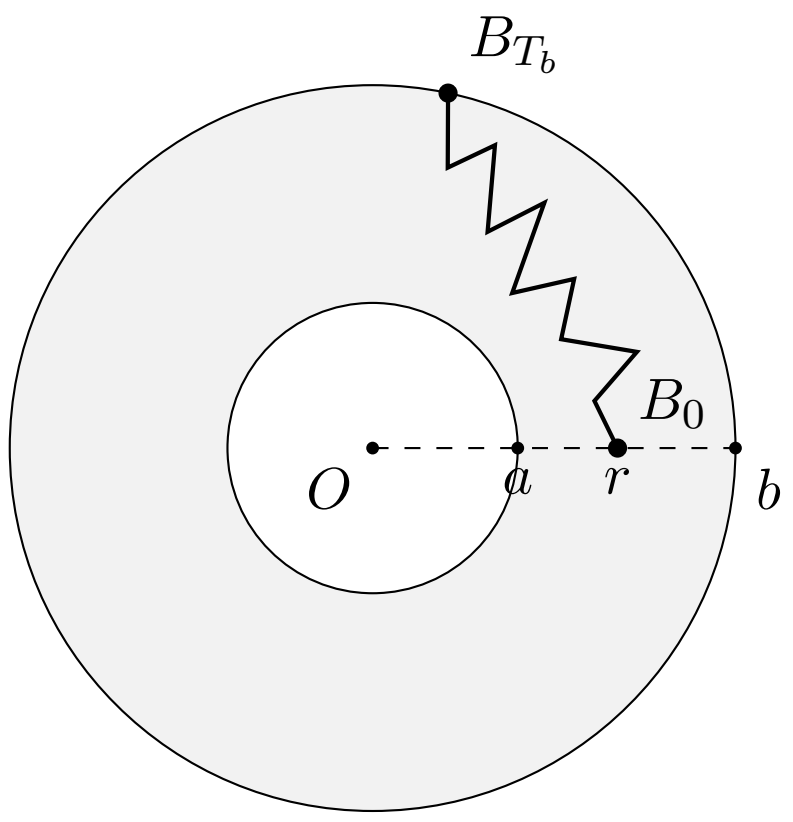
\begin{figure}[htbp]
\centering
\begin{tikzpicture}[scale=1.2]
  \def\a{0.8}
  \def\r{1.35}
  \def\b{2}

  \fill[gray!10,even odd rule]
    (0,0) circle (\b)
    (0,0) circle (\a);
  \draw (0,0) circle (\a);
  \draw (0,0) circle (\b);

  \draw[dashed] (0,0)--(\b,0);
  \fill (0,0) circle (1pt) node[below left] {$O$};
  \fill (\a,0) circle (1pt) node[below] {$a$};
  \fill (\r,0) circle (1.5pt) node[below] {$r$};
  \fill (\b,0) circle (1pt) node[below right] {$b$};
  \node[above right] at (\r,0) {$B_0$};

  \draw[thick]
    (0:\r)--(12:1.25)--(20:1.55)
    --(30:1.20)--(40:1.45)
    --(48:1.15)--(55:1.65)
    --(62:1.35)--(68:1.8)
    --(75:1.6)--(78:\b);

  \fill (78:\b) circle (1.5pt)
    node[above right] {$B_{T_b}$};
\end{tikzpicture}
\caption{Exit of planar Brownian motion from
an annulus.}
\label{fig:annulus}
\end{figure}

For every integer $d\geq2$, formula~\eqref{eq:annulus-exit}
gives the probability that a Brownian motion starting
at distance $r$ from the origin exits an
annulus centered at the origin through its outer sphere. The planar case
is illustrated in Figure~\ref{fig:annulus}. In this case,
the exit probability reduces to $\log(r/a)/\log(b/a)$.

\begin{proposition}[Compatibility of potentials]\label{lem:potentials}
The logarithmic potential $\log r$ is an increasing
convex function of $\Phi_d(r)$ if and only if $d\geq2$.
Equivalently,
\[
 \log r\leq\E_r\log R_{T_a\wedge T_b}
\]
for every $0<a<r<b$ if and only if $d\geq2$.
Equality holds for $d=2$, and the inequality is
strict for $d>2$.
\end{proposition}
\begin{proof}[Proof of Proposition \ref{lem:potentials}]
Writing $y=\Phi_d(r)$, we have
\[
 \frac{\dd}{\dd y}\log r=r^{d-2},\quad
 \frac{\dd^2}{\dd y^2}\log r=(d-2)r^{2d-4}.
\]
The assertion follows from Lemma~\ref{lem:scale}.
\end{proof}

However, the logarithmic potential is singular at the origin.
To apply logarithmic potentials to determinants that may vanish, we consider the regularized logarithms. For $b>0$, direct differentiation gives
\begin{equation}\label{eq:besq-log}
 \mathcal B_d\log(b+s)
 =\frac{db+(d-2)s}{(b+s)^2}.
\end{equation}
Thus $\mathcal B_d\log(b+s)\geq0$ for all $s>0$
if and only if $d\geq2$.

For an integer $d\geq1$, consider the regularized
logarithmic potential $\psi_d$ defined in
\eqref{eq:regularized-potential}.
With $s=c|u-u_0|^2$, the radial identity
\[
 \frac12\Delta_u\!\left[f(c|u-u_0|^2)\right]=c(\mathcal B_d f)(s)
\]
for $f\in C^2([0,\infty))$, together with
\eqref{eq:besq-log}, yields
\begin{equation}\label{eq:euclidean-log}
 \Delta_u\psi_d(u)
 =\frac{2c\{db+(d-2)c|u-u_0|^2\}}
        {(b+c|u-u_0|^2)^2}.
\end{equation}
For $c>0$, the function $\psi_d$ is subharmonic
if and only if $d\geq2$. For $c=0$, it is constant.

The advantage of logarithmic subharmonicity is
that it allows us to take nonnegative linear combinations before exponentiation.

\begin{lemma}\label{lem:log-subharmonicity}
Let $P_t\ (t\geq0)$ be a conservative transition
semigroup and let $g_1,\ldots,g_m$ be strictly positive
measurable functions satisfying
$P_t|\log g_j|<\infty$ and
$P_t\log g_j\geq\log g_j$ for every $t\geq0$.
Then $f=\prod_{j=1}^m g_j^{\nu_j}$ is subharmonic
for all  $\nu_1,\ldots,\nu_m\geq0$.
\end{lemma}

\begin{proof}[Proof of Lemma~\ref{lem:log-subharmonicity}]
Jensen's inequality and order preservation give
\[
 P_t f\geq \exp\Biggl(\sum_{j=1}^m\nu_j P_t\log g_j\Biggr)
 \geq \exp\Biggl(\sum_{j=1}^m\nu_j\log g_j\Biggr)=f. \qedhere
\]
\end{proof}
In applications, logarithmic subharmonicity is often
established by It\^o's formula. Once the drift of $\log g_j(X_t)$ is shown to be
nonnegative, it remains to pass from a local
submartingale to a submartingale. If
\[
 \E_x\sup_{0\leq t\leq T}|\log g_j(X_t)|<\infty
 \quad (T<\infty),
\]
the process $\log g_j(X_t)$ is a submartingale by localization and dominated
convergence. The hypotheses of
Lemma~\ref{lem:log-subharmonicity} then follow.

\subsection{Invariant rays and Schur complements}\label{sec:bessel}

In the continuous shape range, we have seen that
the comparison by rows leads to the noncentrality
matrix $t^2v^*v$.
Let $P_V^{(\alpha)}$ be the Wishart operators
introduced in Appendix~\ref{sec:wishart-semigroup}.
By~\eqref{eq:shifted-row-transform}
and~\eqref{eq:cond-transform},
\begin{equation}\label{eq:cond-one-row}
 \E F_M(Y+(g+tv)^*(g+tv))
 =P_{\Psi/\beta}^{(\alpha)}F_M(t^2v^*v),
\end{equation}
where $\alpha=n-\beta(r-1)/2$. The required comparison is therefore
along the ray generated by $H = v^*v$. 

The operators $P_V^{(\alpha)}$ and $P_H^{(\alpha)}$
commute for every $V,H\in\PP_q$.
To apply Proposition~\ref{cor:ray}, we shall show that the ray
generated by $H$ is invariant under $P_{tH}^{(\alpha)}$
and that $F_M$ is subharmonic with respect to
this semigroup.

\begin{lemma}\label{lem:wishart-ray}
For $H=cc^*\ne0$, the ray $\{xH:x\geq0\}$ is
invariant under the Wishart semigroup $P_{tH}\ (t\geq0)$.
Starting at $xH$, the restriction of the semigroup
to this ray admits a continuous realization $R_tH$,
where $R_t$ is a squared Bessel process of dimension
$2\alpha$ starting at $x$. 
\end{lemma}

\begin{proof}[Proof of Lemma \ref{lem:wishart-ray}]
With $s=c^*Tc$, the determinant and inverse
formulas reduce~\eqref{eq:cond-transform} to
\[
 P_{tH}e_T(xH)
 =(1+2ts)^{-\alpha}
 \exp\left(-\frac{xs}{1+2ts}\right).
\]
This is the Laplace transform of $R_tH$, where
$R_t$ is a squared Bessel process of dimension
$2\alpha$ starting at $x$. Both assertions
follow from uniqueness of Laplace transforms.
\end{proof}

By Lemmas~\ref{lem:wishart-continuity},
\ref{lem:wishart-ray} and~\ref{lem:besq-hitting}, all hypotheses of Proposition~\ref{cor:ray} for $F_M$ have been verified except for subharmonicity with respect to $P_{tH}$.


Suppose, after a permutation of coordinates, that
the lower-right principal block of $H$ vanishes.
For any Hermitian matrix $T_0$ of the corresponding
order, put
\[
 H=\begin{bmatrix}H_{11}&H_{12}\\ H_{12}^*&0\end{bmatrix},
 \qquad
 T=\begin{bmatrix}0&0\\0&T_0\end{bmatrix}.
\]
Then $\tr(TH)=0$ and $THT=0$.
By Lemma~\ref{lem:wishart-characteristics},
$\tr(TU_t)-\tr(TU_0)$ is a continuous martingale
with zero quadratic variation and is therefore
identically zero.
Since $T_0$ is arbitrary, the lower-right
principal block of $U_t$ remains constant. 

This observation leads to the following Schur
complement calculation.

\begin{lemma}
\label{lem:schur-dynamics}
Let $U_t\ (t\geq0)$ be a Wishart diffusion on $\PP_d$
with direction $H=e_1e_1^*$.
For $M\in\PP_d$ and $\varepsilon>0$, partition
\[
 U_t+M+\varepsilon I
 =\begin{bmatrix}
 a_t&b_t^*\\
 b_t&D
 \end{bmatrix},
\]
where $D$ is of order $d-1$ and is constant in time.
Let $C_0$ be the lower-right principal block of
$U_0$, and put
\[
 z_t=a_t-b_t^*D^{-1}b_t,\quad
 w_t=(1,-b_t^*D^{-1})^*.
\]
Then the Schur complement $z_t$ satisfies
\begin{align}
 \dd z_t
 &=\{2\alpha-\beta\tr(D^{-1}C_0)\}\,\dd t+\dd N_t,
 \label{eq:cond-schur-drift}\\
 \dd\langle N\rangle_t
 &=4w_t^*U_tw_t\,\dd t,
 \label{eq:cond-schur-bracket}
\end{align}
where $N_t$ is a continuous local martingale.
Moreover, the drift rate of $\log z_t$ is bounded
below by $\{2\alpha-\beta(d-1)-2\}/{z_t}$.
\end{lemma}

\begin{proof}[Proof of Lemma~\ref{lem:schur-dynamics}]
Since only the first row and column of $U_t$ vary,
its lower-right principal block remains equal to $C_0$.
In particular, $D$ is constant. It\^o's formula gives
\[
 \dd z_t=w_t^*(\dd U_t)w_t
          -(\dd b_t)^*D^{-1}\dd b_t.
\]
The first term has drift $2\alpha\,\dd t$,
since $e_1^*w_t=1$.
By polarization of~\eqref{eq:cond-bracket},
the quadratic variation of its martingale part has differential $4w_t^*U_tw_t\,\dd t$.

To evaluate the second-order term, take
\[
 T=\frac12
 \begin{bmatrix}
 0&v^*\\
 v&0
 \end{bmatrix},
 \quad v\in\F^{d-1}.
\]
Formula~\eqref{eq:cond-bracket} gives
\[
 \dd\langle\operatorname{Re}(v^*b)\rangle_t
 =v^*C_0v\,\dd t.
\]
By polarization, $\dd b_t$ has the same real covariance as $C_0^{1/2}\dd W_t$, where the real components of the $\F$-valued Brownian column $W$ are independent standard Brownian motions.
Consequently,
\[
 (\dd b_t)^*D^{-1}\dd b_t
 =\beta\tr(D^{-1}C_0)\,\dd t.
\]
This establishes~\eqref{eq:cond-schur-drift}
and~\eqref{eq:cond-schur-bracket}.

Since $D\geq C_0$, we have $\tr(D^{-1}C_0)\leq d-1$.
Moreover,
\[
 w_t^*U_tw_t
 \leq w_t^*(U_t+M+\varepsilon I)w_t=z_t.
\]
It\^o's formula for $\log z_t$ gives the drift bound
\[
 \frac{2\alpha-\beta\tr(D^{-1}C_0)}{z_t}
 -\frac{2w_t^*U_tw_t}{z_t^2}
 \geq\frac{2\alpha-\beta(d-1)-2}{z_t}.
\]
For $d=1$, the lower blocks are empty and the trace
term is understood to be zero.
\end{proof}

Since $\det(U_t+M+\varepsilon I)=\det(D)z_t$
and $D$ is constant, the preceding estimate
shows that the regularized logarithmic determinant
is a local submartingale whenever
$2\alpha-\beta(d-1)\geq2$.

If $M=0$, $\varepsilon=0$ and $C_0>0$, then $D=C_0$ and $w_t^*U_tw_t=z_t$.
Then~\eqref{eq:cond-schur-drift}
and~\eqref{eq:cond-schur-bracket} then reduce to
\[
 \begin{aligned}
 \dd z_t &=\{2\alpha-\beta(d-1)\}\,\dd t+\dd N_t,\\
 \dd\langle N\rangle_t &=4z_t\,\dd t.
 \end{aligned}
\]
Hence the Schur complement is itself a squared Bessel
process of dimension $2\alpha-\beta(d-1)$. Thus the sufficient condition above is that this effective dimension be at least two, the threshold in Proposition~\ref{lem:potentials}.

We now apply the preceding Schur complement
calculation to the function $F_M$
in~\eqref{eq:cond-observable}.
Note that $F_M$ is nonnegative, continuous
and of polynomial growth.

\begin{lemma}
\label{lem:cond-ray}
Suppose that $\alpha>0$, $2\alpha\geq\beta(q-1)$
and $2\alpha-\beta(p_j-1)\geq2$
for $j=2,\ldots,m$. Then, for every $V\in\PP_q$
and $c\in\F^q$, $P_V^{(\alpha)}F_M(xcc^*)$
is nondecreasing in $x\ge0$.
\end{lemma}
\begin{proof}[Proof of Lemma~\ref{lem:cond-ray}]
Let $H=cc^*$ with $c\ne0$, fix $\varepsilon>0$,
and let $U_t$ be the Wishart diffusion with
semigroup $P_{tH}$ starting at an arbitrary
$U\in\PP_q$.

By the affine transform~\eqref{eq:cond-transform},
each principal block $U_{jj}(t)$ is again a
Wishart diffusion with shape parameter $\alpha$
and direction $c_jc_j^*$, where $c_j$ is the
corresponding block of $c$.

If $c_j\ne0$, a unitary change of basis sends $c_jc_j^*$ to
$\|c_j\|^2 e_1e_1^*$. Thus time $t$ in the block corresponds to time $\|c_j\|^2t$ in the direction $e_1e_1^*$. By Lemma~\ref{lem:schur-dynamics}, the
process
$\log\det(M_{jj}+U_{jj}(t)+\varepsilon I)$
is a local submartingale. If the direction vanishes
on this block, the process is constant.

The absolute value of each
regularized logarithmic determinant is bounded by a constant
times $1+\tr U_t$. By Lemma~\ref{lem:wishart-characteristics},
$\tr U_t$ is a nonnegative submartingale.
The moment estimate~\eqref{eq:cond-moments} and Doob's
$L^2$ maximal inequality give
\[
 \E_U\sup_{0\leq t\leq T}\tr U_t
 \leq 2\left[\E_U(\tr U_T)^2\right]^{1/2}
 <\infty \ (T<\infty).
\]
Hence, by localization and dominated convergence,
the regularized logarithmic determinants are submartingales.

Lemma~\ref{lem:log-subharmonicity} therefore gives $P_{tH}F_{M+\varepsilon I}
 \geq F_{M+\varepsilon I}.$
A polynomial bound allows us to let $\varepsilon\downarrow0$
by dominated convergence, yielding $P_{tH}F_M\geq F_M$.
The assertion follows from Proposition~\ref{cor:ray}.
The case $c=0$ is immediate.
\end{proof}

\section{Wishart product inequalities}\label{sec:wishart}
We now turn to the direct row comparison
in~\eqref{eq:direct-row-comparison}.
For a single determinant,
Anderson's inequality \cite[Theorem~1]{Anderson1955} gives the required comparison
over all three fields.
For products of determinants, we use logarithmic
subharmonicity in the complex and quaternionic cases.

\begin{lemma}[Anderson's inequality]
\label{lem:anderson}
Let $g$ be a centered Gaussian vector in $\R^d$
and let $K\subset\R^d$ be a convex Borel set symmetric about the origin.
Then, for every $v\in\R^d$, the function $\mathbb P(g+tv\in K)$
is nonincreasing in $ t \ge 0$.
\end{lemma}

\begin{lemma}\label{lem:single-det}
Let $D\in\PP_d$, let $g$ be a centered
$\F$-Gaussian row and let $v$ be a deterministic row.
Then the function $\det\bigl(D+(g+tv)^*(g+tv)\bigr)$
is stochastically nondecreasing in $t\geq0$.
\end{lemma} 
\begin{proof}[Proof of Lemma~\ref{lem:single-det}]
Suppose first that $D>0$ and that $g$ has a
positive definite covariance matrix. By the matrix determinant lemma,
\[
 \det(D+(g+tv)^*(g+tv))
 =\det D\{1+(g+tv)D^{-1}(g+tv)^*\}.
\]
For every $c\geq0$, the ellipsoid $\{u\in\F^{1\times d}:uD^{-1}u^*\leq c\}$
is convex and symmetric about the origin.
By Anderson's inequality, the probability
\[
 \mathbb P((g+tv)D^{-1}(g+tv)^*\leq c)
\]
is nonincreasing in $t\geq0$.
This proves the stochastic monotonicity for $D>0$.

For the general case, let $\xi$ be an independent
standard $\F$-Gaussian row and replace $g$ by
$g+\sqrt{\varepsilon}\xi$ and $D$ by
$D+\varepsilon I_d$.
Since stochastic order is preserved under weak
convergence, letting $\varepsilon\downarrow0$
gives the assertion.
\end{proof}

\begin{lemma}\label{lem:cond-row}
Suppose that $\beta=2$ or $4$, and let $M\in\PP_q$. Let $g\in\F^{1\times q}$ be a random row such that $ag\overset{\mathrm d}{=}g$ for every $a\in\F$ with $|a|=1$. If $\E F_{M+I}(g^*g)<\infty$, then the function $\E F_M((g+tv)^*(g+tv))$ is nondecreasing in $t\ge0$ for every deterministic row $v$.
\end{lemma}

\begin{proof}[Proof of Lemma~\ref{lem:cond-row}]
For $0<\varepsilon\leq1$  and $u\in\F$, put
\[
 G_\varepsilon(u)
 =\E F_{M+\varepsilon I}((g+uv)^*(g+uv)).
\]
For block $j$, put $A=M_{jj}+\varepsilon I$. By the matrix determinant lemma, the regularized determinant is
\[
(\det A)\bigl\{1+\|(g_j+uv_j)A^{-1/2}\|^2\bigr\}.
\]
Put $c=v_jA^{-1}v_j^*$. If $c=0$, the factor in braces is the constant
$1+g_jA^{-1}g_j^*$. If $c>0$, completing the square gives $b+c|u-u_0|^2$, where
\[
 u_0=-\frac{g_jA^{-1}v_j^*}{c},\qquad
 b=1+g_jA^{-1}g_j^*-c|u_0|^2 
\]
The Cauchy--Schwarz inequality gives $b\geq1$. Thus, by~\eqref{eq:euclidean-log}, the logarithm of the factor in braces is subharmonic on $\F\simeq\R^\beta$.

Taking the exponential of the weighted sum of these subharmonic logarithms preserves subharmonicity, as does integration with respect to $g$. Indeed, for $u$ in a bounded set, the matrix $M+\varepsilon I+(g+uv)^*(g+uv)$ is bounded above
in the Loewner order by a constant multiple of
$M+I+g^*g$.
The integrability assumption therefore provides a common integrable
bound on the regularized products and gives continuity of $G_\varepsilon$.

Moreover, the invariance of $g$ under left
multiplication by unit scalars implies that
$G_\varepsilon$ is radial.
For Brownian motion $B_t$ starting at $u$ with $|u|<r$,
let $\tau_r$ be its first hitting time of the
sphere of radius $r$. By Lemma~\ref{lem:besq-hitting},
$\tau_r$ is almost surely finite.
Since $G_\varepsilon$ is bounded on the closed
ball of radius $r$, optional sampling gives
\[
 G_\varepsilon(u)
 \leq \E_u G_\varepsilon(B_{\tau_r})
 =G_\varepsilon(r).
\]
Thus $G_\varepsilon$ is nondecreasing in $|u|$.
Dominated convergence as $\varepsilon\downarrow0$ completes the proof.
\end{proof}

For the conditional moment $G$, retain the assumptions and notation of Section~\ref{sec:conditional}.

\begin{lemma}
\label{prop:cond-moments}
Suppose either that $2n\geq\beta(p-1)$ or that
$\beta\in\{2,4\}$ and $n=k\beta/2$ for an integer
$k\geq\max_j p_j$.
Then $G$ has a continuous version on the positive definite cone that is nondecreasing in the Loewner order.
\end{lemma}

\begin{proof}[Proof of Lemma~\ref{prop:cond-moments}]
Both shape conditions satisfy the hypotheses of Lemma~\ref{lem:cond-regression}. At the endpoint $2n=\beta(p-1)$, the discrete representation applies with $k=p-1\geq r$.
The continuity of $G$ follows from
\eqref{eq:conditional-reduction}, as established in
Section~\ref{sec:conditional}.
By Lemma~\ref{lem:row-reduction}, it remains to
establish the monotonicity in~\eqref{eq:row-comparison}
for every $M\in\PP_q$ and every deterministic row $v$.

If $m=2$, raising the determinant in Lemma~\ref{lem:single-det} to the positive power $\nu_2$ gives the monotonicity in~\eqref{eq:direct-row-comparison}.

If $m\geq3$ and $2n\geq\beta(p-1)$, put
$\alpha=n-\beta(r-1)/2$.
Since $2\alpha\geq\beta q$ and $p_j\leq q-1$
for $j=2,\ldots,m$, we have
\[
 2\alpha-\beta(p_j-1)
 \geq\beta(q-p_j+1)\geq2.
\]
Lemma~\ref{lem:cond-ray} therefore applies. Together
with~\eqref{eq:cond-one-row}, it gives the required monotonicity in~\eqref{eq:row-comparison}.

In the remaining case, $\beta\in\{2,4\}$ and
$n=k\beta/2$.
The centered Gaussian row $g$ is invariant in law
under left multiplication by elements of unit modulus.
Moreover, the polynomial growth of $F_{M+I}$
and the Gaussian moment bounds ensure that
$\E F_{M+I}(g^*g)<\infty$.
Lemma~\ref{lem:cond-row} therefore gives the
monotonicity in~\eqref{eq:direct-row-comparison}.

In the two cases that use~\eqref{eq:direct-row-comparison}, conditioning on $Y$ gives~\eqref{eq:row-comparison} as well. Lemma~\ref{lem:row-reduction} now completes the proof.
\end{proof}

This verifies the sufficient condition for separation in Lemma~\ref{lem:spectral-reduction}. To characterize equality, we determine when $G$ is constant.
\begin{lemma}
\label{lem:cond-nonconstancy}
Under the assumptions of
Lemma~\ref{prop:cond-moments}, suppose further
that $\Sigma_{jj}>0$ for $j=2,\ldots,m$.
If $\Sigma_{AB}=0$, then $G$ is constant.
Otherwise, $G(t\Sigma_{AA})\to\infty$ as $t\to\infty$.
\end{lemma}

\begin{proof}[Proof of Lemma \ref{lem:cond-nonconstancy}]
If $\Sigma_{AB}=0$, then $C=0$ and $G$ is constant by Lemma~\ref{lem:cond-regression}. Suppose now that $\Sigma_{AB}\ne0$ and put $D=\Sigma_{AA}$. For $t\geq1$, vary the covariance of the first block while keeping $C$ and $\Psi$ fixed by setting
\[
 \Sigma(t)=
 \begin{bmatrix}
 tD&tDC\\
 tC^*D&\Psi+tC^*DC
 \end{bmatrix}.
\]
Since its Schur complement with respect to $tD$
is $\Psi\geq0$, we have $\Sigma(t)\geq0$.
Let $S(t)\sim\W_p^{(\beta)}(n,\Sigma(t))$.

The first block $S_{11}(t)$ has the same law as $tA$. The conditional law of the complementary block uses the same $C$ and $\Psi$ for every $t$.
Lemma~\ref{lem:cond-regression} therefore gives
\[
 \E G(tA)
 =\E\prod_{j=2}^m\det(S_{jj}(t))^{\nu_j}.
\]
Under either shape condition, we have
$2n>\beta(p_j-1)$ for $j=2,\ldots,m$.
Thus each $\log\det S_{jj}(t)$ has finite moments. Jensen's inequality gives
\[
 \begin{aligned}
 \log\E G(tA)
 &\geq\sum_{j=2}^m\nu_j\E\log\det S_{jj}(t)\\
 &=c+\sum_{j=2}^m\nu_j\log\det\Sigma_{jj}(t),
 \end{aligned}
\]
where $c$ is independent of $t$ by congruence of the marginal Wishart laws.

Let $C_j$ denote the submatrix of $C$ corresponding
to block $j$. Since
$\Sigma_{jj}(t)=\Sigma_{jj}+(t-1)C_j^*DC_j$,
its determinant is nondecreasing in $t$ and
tends to infinity whenever $C_j\ne0$.
Hence $\E G(tA)\to\infty$ if $C\ne0$.

Every $a>0$ satisfies $a\leq tD$ for some $t>0$. Thus, if $G(tD)$ were bounded on this ray, Loewner monotonicity would bound $G$ on the whole positive definite cone, contradicting the preceding limit. 
\end{proof}

We now combine the spectral reduction in
Section~\ref{sec:association} with the conditional
monotonicity established above.

\begin{lemma}
\label{lem:spectral}
Suppose that $m\geq2$, $\Sigma\geq0$, and $\Sigma_{jj}>0$ and $\nu_j>0$ for every $j=1,\ldots,m$. Assume either that $2n\geq\beta(p-1)$ or that $\beta\in\{2,4\}$ and $n=k\beta/2$ for an integer $k\geq\max_j p_j$.
Then the separation inequality~\eqref{eq:sep}
holds, with equality if and only if
$\Sigma_{1j}=0$ for every $j=2,\ldots,m$.
\end{lemma}

\begin{proof}[Proof of Lemma~\ref{lem:spectral}]
Put $A=S_{11}$ and $D=\Sigma_{11}$.
By Lemma~\ref{prop:cond-moments}, the conditional
moment $G$ defined in~\eqref{eq:conditional-moment}
admits a continuous, Loewner-nondecreasing version. Under either shape condition, $2n>\beta(p_1-1)$. Lemma~\ref{lem:spectral-reduction} therefore gives~\eqref{eq:sep},
with equality if and only if $G$ is constant.
By Lemma~\ref{lem:cond-nonconstancy}, this is equivalent
to $\Sigma_{AB}=0$.
\end{proof}

The Wishart Laplace transform shows that every principal submatrix has a Wishart distribution with the same shape parameter and the corresponding principal submatrix of $\Sigma$ as its covariance parameter. We may therefore omit all blocks with $\nu_j=0$. If at most one block remains, the assertions are immediate.

\begin{proof}[Proof of Theorem~\ref{thm:real}]
 If some retained $\Sigma_{jj}$ is singular, then $S_{jj}$ is
 singular almost surely and both sides of~\eqref{eq:wpi} vanish.
 Otherwise all retained $\Sigma_{jj}$ are positive definite.
 Successive applications of Lemma~\ref{lem:spectral} to the complementary principal marginals give~\eqref{eq:wpi}.

Assume that $\Sigma_{ij}\ne0$ for some $i\ne j$.
Relabel the blocks so that $\Sigma_{1j}\ne0$
for some $j\geq2$.
By Lemma~\ref{lem:spectral}, the first
separation is strict.
Since all marginal determinant moments are positive, inequality~\eqref{eq:wpi} is strict.

Conversely, if $\Sigma_{ij}=0$ for all $i\ne j$
with $\nu_i\nu_j\ne0$, the joint Laplace transform
of the corresponding diagonal blocks factors
into their marginal transforms.
These blocks are therefore independent,
giving equality.
\end{proof}

\begin{proof}[Proof of Corollary~\ref{cor:ch}]
If some $\Sigma_{jj}$ is singular, then $S_{jj}$
is singular almost surely by the congruence
property of the Wishart distribution.

At a discrete shape $n=k\beta/2$, write $S=X^*X$,
where $X$ has $k$ independent centered Gaussian rows.
Let $X_j$ denote the $j$th block of columns of $X$.
Since $S_{jj}=X_j^*X_j$, we have
$\rank(S_{jj})\leq k$, so $S_{jj}$ is singular whenever $k<p_j$. In either case, both sides of~\eqref{eq:wpi} vanish.

Otherwise, all retained $S_{jj}$ are positive definite almost surely. Successive applications of Lemma~\ref{lem:spectral}
give inequality~\eqref{eq:wpi}.
The equality condition follows as in the proof
of Theorem~\ref{thm:real}.
\end{proof}

\begin{proof}[Proof of Corollary~\ref{thm:gauss}]
Take $S=Z^*Z$ and $n=2$. Since $S_{jj}=|Z_j|^2$, applying Corollary~\ref{cor:ch} to scalar blocks with exponents $\nu_j/2$ gives~\eqref{eq:gpi}. If $\Sigma_{jj}>0$ for every $j$, then all diagonal entries of $S$ are positive almost surely, and the equality condition follows from Corollary~\ref{cor:ch}.
\end{proof}

\section{Further product inequalities}
\label{sec:non-gaussian}

In this section, we extend the preceding results in
three directions. We first establish the Wishart product
inequality at the remaining admissible shapes in the
real case. We then replace the determinant powers by
a more general class of functions and extend the
comparison by rows to certain non-Gaussian
random matrices.

\subsection{Singular real Wishart distributions}
\label{sec:singular-real}

At discrete shapes below the continuous range, the proof over
$\C$ and $\HH$ uses the comparison by rows in Lemma~\ref{lem:cond-row}.
In real dimension one, logarithmic subharmonicity fails, and the
example in the Discussion shows that conditional monotonicity can fail
as well. To prove the product inequality at these real shapes, we
use the strong Gaussian product inequality of Ouimet and
Greaves~\cite{OuimetGreaves2026}, which gives~\eqref{eq:gpi} over
$\mathbb R$. We first derive a consequence for Gaussian quadratic forms.

\begin{lemma}\label{lem:gaussian-quadratic-products}
Let $g=(g_1,\ldots,g_m)$ be a centered $\F$-Gaussian
row partitioned into blocks. For positive semidefinite
matrices $Q_j$ of the corresponding orders,
$b_j\geq0$ and $\nu_j>0$, we have
\[
 \E\prod_{j=1}^m(b_j+g_jQ_jg_j^*)^{\nu_j}
 \geq
 \prod_{j=1}^m\E(b_j+g_jQ_jg_j^*)^{\nu_j}.
\]
\end{lemma}

\begin{proof}[Proof of Lemma~\ref{lem:gaussian-quadratic-products}]
Let $h_j$ be standard $\F$-Gaussian columns of the
corresponding orders and let $\xi_j$ be standard
$\F$-Gaussian variables, mutually independent
and independent of $g$. Put
\[
 Y_j=g_jQ_j^{1/2}h_j+\sqrt{b_j}\,\xi_j.
\]
Conditionally on the columns $h_1,\ldots,h_m$,
the vector $(Y_1,\ldots,Y_m)$ is centered
$\F$-Gaussian, and the conditional marginal
distribution of $Y_j$ depends only on $h_j$.
Applying the Gaussian product inequality over the corresponding field with exponents $2\nu_j$
and using the independence of the columns gives
\[
 \begin{aligned}
 \E\prod_{j=1}^m|Y_j|^{2\nu_j}
 &\geq\ E\prod_{j=1}^m
       \E[|Y_j|^{2\nu_j}\mid h_j] =\prod_{j=1}^m\E|Y_j|^{2\nu_j}.
 \end{aligned}
\]

On the other hand, conditionally on $g$, the variables
$Y_j$ are independent and each has the distribution
of $(b_j+g_jQ_jg_j^*)^{1/2}\xi_j$.
Putting $c_j=\E|\xi_j|^{2\nu_j}$, we obtain $\E|Y_j|^{2\nu_j}
 =c_j\E(b_j+g_jQ_jg_j^*)^{\nu_j}.$
Moreover, the conditional independence of the
variables $Y_j$ gives
\[
 \E\prod_{j=1}^m|Y_j|^{2\nu_j}
 =\Bigl(\prod_{j=1}^m c_j\Bigr)
  \E\prod_{j=1}^m(b_j+g_jQ_jg_j^*)^{\nu_j}.
\]
Substituting these identities into the preceding
inequality and cancelling the common factors
gives the assertion.
\end{proof}

\begin{theorem}[WPI in the real discrete shape range]\label{prop:singular-real}
Let $S\sim\W_p^{(1)}(k/2,\Sigma)$, where $\Sigma\geq0$
and $k=0,\ldots,p-1$.
Then inequality~\eqref{eq:wpi} holds.
\end{theorem}

\begin{proof}[Proof of Theorem~\ref{prop:singular-real}]
After discarding the blocks with zero exponents,
the assertion is immediate if at most one block
remains.
Write $S=X^*X$, where $X$ has $k$ independent
centered Gaussian rows of covariance $\Sigma$.
If $k<\max_j p_j$, at least one retained principal
block is singular almost surely, so both sides
of~\eqref{eq:wpi} vanish.
We may therefore assume that $k\geq\max_j p_j$.

Fix $i$ and condition on all rows except the $i$th row $g_i$.
Write $g_{\ell j}$ for block $j$ of row $g_\ell$. The contribution of
the other rows to the $j$th principal block is $B_j=\sum_{\substack{1\leq\ell\leq k\\ \ell\ne i}}g_{\ell j}^*g_{\ell j}.$
The determinant identity
\[
 \det(B_j+g_{ij}^*g_{ij})
 =\det B_j+g_{ij}\operatorname{adj}(B_j)g_{ij}^*
\]
and the positive semidefiniteness of $\operatorname{adj}(B_j)$ allow us to apply Lemma~\ref{lem:gaussian-quadratic-products} with $b_j=\det B_j$ and $Q_j=\operatorname{adj}(B_j)$.

Thus replacing the blocks of $g_i$ by independent Gaussian vectors with the same marginal laws, also independent of the other rows, can only decrease the conditional product moment. After integrating over the other rows, we may repeat the replacement for each row in turn. The resulting principal blocks are independent and retain their original marginal laws. This gives~\eqref{eq:wpi}.
\end{proof}

\subsection{Functional product inequalities}

The preceding product inequalities involve positive
powers of principal determinants. Convex, nondecreasing transformations of logarithmic potentials preserve subharmonicity. This allows us to replace the determinant powers by a broader class of functions of principal minors. Spectral association handles the first factor in the separation, so that factor requires only monotonicity under the standing continuity and growth assumptions below.

For $S\in\PP_p$, write $D_j(S)=\det S_{jj}$.
In this and the next subsection, let
$f_j:[0,\infty)\to[0,\infty)$ be continuous,
nondecreasing functions of polynomial growth
and put $f(S)=\prod_{j=1}^m f_j(D_j(S))$.
The additional condition used in the comparisons is
\begin{equation}\label{eq:functional-condition}
 f_j(x)>0\ (x>0),\quad
 s\mapsto\log f_j(e^s)\text{ is convex on }\R.
\end{equation}

\begin{theorem}[Functional product inequality]
\label{prop:functional}
Let $S\sim\W_p^{(\beta)}(n,\Sigma)$ with $\Sigma\geq0$,
where $n\geq(p-1)/2$ for $\beta=1$ and
$n\in\beta\mathscr A_p$ for $\beta\in\{2,4\}$.
If all but at most one of the functions $f_j$
satisfy~\eqref{eq:functional-condition}, then
\begin{equation}\label{eq:functional}
 \E\prod_{j=1}^m f_j(D_j(S))
 \geq\prod_{j=1}^m\E f_j(D_j(S)).
\end{equation}
\end{theorem}

\begin{proof}[Proof of Theorem~\ref{prop:functional}]
First suppose that $\Sigma>0$. At a discrete shape
$n=k\beta/2$, each block of order $p_j>k$ contributes
the constant $f_j(0)$ on both sides. If any such constant is zero, the inequality follows immediately. Otherwise cancel these factors  and pass to the principal marginal on the retained blocks. The two sides agree if at most one block remains.

Place the exceptional function first, if there is one,
and put $\phi_j(s)=\log f_j(e^s)$ for $j\geq2$.
Each $\phi_j$ is convex and nondecreasing.
In the argument of Lemma~\ref{lem:cond-ray}, let $L_j(t)$ be
the logarithm of the regularized determinant in block $j$.
For fixed $\varepsilon>0$, this determinant is at least $\varepsilon^{p_j}$. Polynomial growth and~\eqref{eq:cond-moments} therefore ensure that $\phi_j(L_j(t))$ is integrable. For $0\leq s\leq t$, conditional Jensen's inequality and the submartingale property of $L_j$ give
\[
 \begin{aligned}
 \E[\phi_j(L_j(t))\mid\mathscr F_s]
 &\geq\phi_j(\E[L_j(t)\mid\mathscr F_s])\geq\phi_j(L_j(s)) .
 \end{aligned}
\]
Thus each $\phi_j(L_j)$ is a submartingale. Applying conditional Jensen's inequality to the exponential of their sum shows that $\exp\{\sum_{j=2}^m\phi_j(L_j)\}$ is also a submartingale.

The same convex, nondecreasing composition preserves
the subharmonicity used in the proof of
Lemma~\ref{lem:cond-row}.
In both comparisons, polynomial domination allows
us to remove the regularization.
When only one block remains to be compared, the
stochastic comparison in Lemma~\ref{lem:single-det}
applies directly.

The argument of Lemma~\ref{prop:cond-moments} now gives a continuous, Loewner-nondecreasing version $G_f$ of the conditional expectation of
$\prod_{j=2}^m f_j(D_j(S))$ given $S_{11}=a$. By the observation following
Lemma~\ref{lem:spectral-reduction}, we obtain
\[
 \Cov(f_1(\det A),G_f(A))\geq0.
\]
This gives the required separation inequality.
Iteration on the complementary principal marginals
proves~\eqref{eq:functional}.

For singular $\Sigma$, replace $\Sigma$ by
$\Sigma_N=\Sigma+N^{-1}I_p$ and let $N\to\infty$.
The continuity and polynomial growth of the
functions $f_j$, together with the Wishart moment
bounds, give convergence of both sides
of~\eqref{eq:functional}.
\end{proof}

For example, the function $f_j(x)=(a_j+x^{b_j})^{c_j}$
satisfies~\eqref{eq:functional-condition} for $a_j\geq0$
and $b_j,c_j>0$. The exceptional function may be bounded,
as in $f_1(x)=x/(1+x)$.

\subsection{Non-Gaussian extensions}

The comparison by rows in Lemma~\ref{lem:cond-row}
extends to the functions satisfying~\eqref{eq:functional-condition}
by the preceding argument. This suggests that the comparison of product moments
can be extended beyond Gaussian matrices.
In particular, suitable invariance properties of
the joint distribution of the rows allow us to
carry out successive row comparisons without
assuming independence. We first establish
monotonicity with respect to both dilations
and deterministic shifts.

A random $k\times p$ matrix $X$ over $\F$ is said
to have diagonal unitary invariance if
$QX\overset{\mathrm d}{=}X$ for every matrix
$Q=\diag(u_1,\ldots,u_k)$ with $u_i\in\F$
and $|u_i|=1$ for $i=1,\ldots,k$.
It is said to be left unitarily invariant if
the same identity holds for every $k\times k$
matrix $Q$ satisfying $Q^*Q=I_k$.

\begin{proposition}[Monotonicity under unitary invariance]
\label{prop:noncentral-monotonicity}
Suppose that $\beta=2$ or $4$ and that every $f_j$
satisfies~\eqref{eq:functional-condition}.
Let $X$ be a random $k\times p$ matrix over $\F$ such that
$\E f(I_p+X^*X)<\infty$.
For $s\geq0$ and deterministic $M\in\F^{k\times p}$, put $h_s(M)=\E f((sX+M)^*(sX+M))$.

If $X$ has diagonal unitary invariance, then $s\mapsto h_s(M)$ is nondecreasing for each $M$. Moreover, for each $s\geq0$ and each
deterministic $k\times p$ matrix $V$, the function
\[
 (t_1,\ldots,t_k)\mapsto h_s(\diag(t_1,\ldots,t_k)V)
\]
is coordinatewise nondecreasing on $[0,\infty)^k$.

If $X$ is left unitarily invariant, then
\begin{equation}\label{eq:noncentral-monotonicity}
 M_0^*M_0\leq M_1^*M_1
 \quad\Longrightarrow\quad h_s(M_0)\leq h_s(M_1).
\end{equation}
\end{proposition}

\begin{proof}[Proof of Proposition~
\ref{prop:noncentral-monotonicity}]
Since the \(f_j\) have polynomial growth, choose $b_j\geq0$ such that
$\log f_j(e^u)\leq b_ju+O(1)$ as $u\to\infty$.
By convexity, the slopes of $\log f_j(e^u)$
are bounded above by $b_j$. It follows that
$f_j(cx)\leq c^{b_j}f_j(x)$ for $c\geq1$ and $x>0$.
The inequality extends to $x=0$ by continuity.
Consequently,
$f(cS)\leq c^af(S)$ for $c\geq1$, where $a=\sum_{j=1}^m p_jb_j$.

For $s$ and $M$ in bounded sets, we have
\[
 (sX+M)^*(sX+M)
 \leq2s^2X^*X+2M^*M
 \leq c(I_p+X^*X)
\]
for some constant $c\geq1$.
Since $f$ is nondecreasing in the Loewner order,
the integrand defining $h_s(M)$ is bounded by
$c^af(I_p+X^*X)$.
The integrability assumption and dominated convergence
therefore imply that $h_s(M)$ is finite and jointly
continuous in $(s,M)$.

To prove the first assertion, allow the factors of the rows of $X$ and $V$ to vary separately. Write $x_i$ and $v_i$ for their $i$th rows. Let $u$ be independent of $X$ and uniformly distributed on the unit scalars in $\F$. By diagonal unitary invariance, replacing $x_i$ by $ux_i$ preserves the joint law of the rows of $X$.

Condition on $X$ and let $B$ be the contribution of all rows except the $i$th, with their multipliers held fixed. For every row $w$ and unit scalar $u$, we have $(\bar u w)^*(\bar u w)=w^*w$. Hence
\[
 \begin{aligned}
 &\E_u f(B+(su x_i+tv_i)^*(su x_i+tv_i))\\
 &\qquad=\E_u f(B+(sx_i+t\bar u v_i)^*
                         (sx_i+t\bar u v_i)).
 \end{aligned}
\]
By the extension of Lemma~\ref{lem:cond-row}
to functions satisfying~\eqref{eq:functional-condition},
the first expression is nondecreasing in $t$,
while the second is nondecreasing in $s$.
Taking expectations over $X$ and varying the
factors of the rows successively proves
the first assertion.

Suppose now that $X$ is left unitarily invariant.
Then $h_s(QM)=h_s(M)$ for every unitary matrix $Q$.
If $M_0^*M_0\leq M_1^*M_1$, we have
$\|M_0x\|\leq\|M_1x\|$ for every $x\in\F^p$.
Hence the map $M_1x\mapsto M_0x$ is well defined
and contractive on the range of $M_1$.
Extending it by zero on the orthogonal complement
of this range gives a contraction $K$ satisfying
$M_0=KM_1$.

By the singular value decomposition, write
$K=Q_0\diag(t_1,\ldots,t_k)Q_1^*$, where $Q_0$
and $Q_1$ are unitary and $0\leq t_i\leq1$.
By unitary invariance and the first assertion,
\[
 \begin{aligned}
 h_s(M_0)
 &=h_s(\diag(t_1,\ldots,t_k)Q_1^*M_1)\\
 &\leq h_s(Q_1^*M_1)=h_s(M_1).
 \end{aligned}
\]
This proves~\eqref{eq:noncentral-monotonicity}.
\end{proof}

\begin{example}[Matrix $t$ distribution]
\label{ex:vector-student}
Suppose that $\beta=2$ or $4$ and $k\geq\max_j p_j$.
Let $G$ be a $k\times p$ Gaussian matrix with independent
centered rows of covariance $\Sigma>0$, and let
$H\sim\chi_\eta^2$ be independent of $G$.
Put $T=\eta/H$ and consider
\[
 S(M)=T(G+M)^*(G+M),\quad S=S(0).
\]
Choose $b_j\geq0$ such that $f_j(x)\leq C_j(1+x)^{b_j}$
and put $a=\sum_{j=1}^m p_jb_j$.

\begin{enumerate}
\item[(i)] If $\eta>2a$ and every $f_j$ satisfies
\eqref{eq:functional-condition}, then $\E f(S(M))$
is nondecreasing in $M^*M$ with respect to the
Loewner order.

\item[(ii)] If $\eta>2a$ and the functions $f_j$,
$j=2,\ldots,m$, satisfy~\eqref{eq:functional-condition},
then the conditional expectation of
$\prod_{j=2}^m f_j(D_j(S))$ given $S_{11}=x$
admits a continuous Loewner-nondecreasing version.
In particular, $S$ satisfies~\eqref{eq:functional}.

\item[(iii)] If $f_j(x)=x^{\nu_j}$ with $\nu_j>0$,
then, for every deterministic matrix $M$,
\begin{equation}\label{eq:student-moments}
 \E f(S(M))<\infty
 \quad\text{if and only if}\quad
 \eta>2\sum_{j=1}^m p_j\nu_j.
\end{equation}
In this range, $S$ satisfies~\eqref{eq:wpi},
strictly if $m\geq2$.
\end{enumerate}

To prove (i), note that $\eta>2a$ ensures
$\E T^a<\infty$ and hence the finiteness of
$\E f(S(M))$ by polynomial growth.
Conditionally on $T=t$, apply
Proposition~\ref{prop:noncentral-monotonicity}
with $X=G$, $s=\sqrt t$ and deterministic matrix
$\sqrt tM$. Integration over $T$ gives the assertion.

For (ii), partition $S$ and $\Sigma$ as in
Section~\ref{sec:conditional}, with $A=S_{11}$,
$D=\Sigma_{AA}$ and $r=p_1$.
Conditionally on $T=t$, we have
$A\sim\W_r^{(\beta)}(k\beta/2,tD)$.
The density of $A$ at $x>0$ depends on $t$ through
the factor $t^{-\beta kr/2}
\exp\{-\beta\tr(D^{-1}x)/(2t)\}$.
Since $T=\eta/H$, Bayes' formula gives
\[
 T\mid A=x\ \overset{\mathrm d}{=}\ T_x
 =\frac{\eta+\beta\tr(D^{-1}x)}{H_1},
 \quad H_1\sim\chi_{\eta+\beta kr}^2.
\]

Partition the Gaussian matrix as $G=(G_A,G_B)$.
Gaussian regression gives $G_B=G_AC+E$, where
$E$ is independent of $(G_A,T)$ and has $k$
independent centered Gaussian rows of covariance
$\Psi$.

Given $(G_A,T)$, choose a unitary matrix $Q$
such that
\[
 Q\sqrt T G_A=\begin{bmatrix}A^{1/2}\\0_{(k-r)\times r}\end{bmatrix}.
\]
By left unitary invariance, $QE$ has the same
distribution as $E$ and is independent of $(G_A,T)$.
Since $B=T G_B^*G_B$, the conditional distribution
of $B$ given $A=x$ is that of $Z_x^*Z_x$, where
\[
 Z_x=\sqrt{T_x}\,E+M_x,\quad
 M_x=\begin{bmatrix}x^{1/2}C\\0_{(k-r)\times(p-r)}\end{bmatrix}.
\]
Here $E$ is taken independent of $H_1$.

For $0<x\leq y$, coupling \(T_x\) and \(T_y\) through the same \(H_1\) gives
$T_x\leq T_y$, while
$M_x^*M_x=C^*xC\leq C^*yC=M_y^*M_y$.
With $X=E$ and $f$ restricted to the remaining
blocks, Proposition~\ref{prop:noncentral-monotonicity}
gives
\[
 h_{\sqrt{T_x}}(M_x)
 \leq h_{\sqrt{T_y}}(M_x)
 \leq h_{\sqrt{T_y}}(M_y).
\]
Integration over $H_1$ gives the required conditional
monotonicity. On compact subsets of the positive
definite cone, $M_x$ is bounded and $T_x$ is bounded
by a constant times $H_1^{-1}$.
Since $\E H_1^{-a}<\infty$, continuity follows
by dominated convergence and Gaussian moment bounds.

The whitened matrix $D^{-1/2}AD^{-1/2}$ is
unitarily invariant. Its ordered eigenvalues are
associated, being products of associated Wishart
eigenvalues and a common independent positive
random variable.
The observation following
Lemma~\ref{lem:spectral-reduction} therefore applies.
Iteration on the complementary principal marginals
proves~\eqref{eq:functional}.

For (iii), put $a=\sum_{j=1}^m p_j\nu_j$.
By independence and homogeneity,
\[
 \E f(S(M))=\E T^a\,\E f((G+M)^*(G+M)).
\]
Since $k\geq\max_j p_j$ and $\Sigma>0$, the Gaussian
moment is finite and strictly positive.
Since $T=\eta/H$ with $H\sim\chi_\eta^2$,
\[
 \E T^a=
 \begin{cases}
 (\eta/2)^a\,\Gamma(\eta/2-a)/\Gamma(\eta/2), & 0\leq a<\eta/2,\\
 +\infty, & a\geq\eta/2.
 \end{cases}
\]
This gives~\eqref{eq:student-moments}.

For $m\geq2$, the nondegeneracy of $T$ and the
positivity of the exponents imply
$\E T^a>\prod_j\E T^{p_j\nu_j}$.
Together with the Wishart product inequality,
this gives
\[
 \begin{aligned}
 \E f(S)&=\E T^a\,\E f(G^*G)\geq \E T^a\prod_{j=1}^m\E D_j(G^*G)^{\nu_j}\\
 & > \Biggl(\prod_{j=1}^m\E T^{p_j\nu_j}\Biggr)\!\Biggl(\prod_{j=1}^m\E D_j(G^*G)^{\nu_j}\Biggr)\\
 & = \prod_{j=1}^m\E D_j(S)^{\nu_j}.
 \end{aligned}
\]
This proves the product inequality and its strictness.
The case $m=1$ is immediate.
\end{example}

\section*{Discussion}

Although the product inequality holds for every
admissible shape parameter, conditional
monotonicity may fail at discrete real shapes.
We conclude with three questions concerning
the limitations of our method.

At a discrete real shape $n=k/2$, we have $S=X^*X$,
where $X$ has $k$ independent centered Gaussian rows
of covariance $\Sigma$. After discarding the blocks
with zero exponents, the inequality is trivial if
$p_j>k$ for some $j$, since the corresponding determinant
vanishes almost surely. Thus the remaining range not
covered by Theorem~\ref{thm:real} is
$\max_j p_j\leq k<p-1$.

In this range, conditional monotonicity may fail.
Let $Z$ be a standard normal variable independent of
a centered Gaussian pair $(G_1,G_2)$ with unit variances
and correlation $3/4$. Put $X_1=Z$, $X_2=G_1+Z/4$
and $X_3=G_2-Z/4$. Wick's formula gives
\[
 \E[X_2^2X_3^2\mid X_1^2=x]
 =\frac{17}{8}-\frac{x}{16}+\frac{x^2}{256}.
\]
This conditional moment is strictly decreasing on
$(0,8)$. Moreover,
\[
 \Cov(X_1^2,X_2^2X_3^2)=-\frac{5}{64}<0,
\]
whereas
\[
 \E\prod_{j=1}^3 X_j^2
 -\prod_{j=1}^3\E X_j^2=\frac{111}{128}>0.
\]
Thus the product inequality may hold even when
both conditional monotonicity and the separation
inequality~\eqref{eq:sep} fail.

This example motivates the following questions.
\begin{question}
For prescribed block sizes $p_2,\ldots,p_m$,
determine all shape parameters $\alpha$ for which
the Wishart kernels $P_V^{(\alpha)}$ exist on
$\PP_q$ for every $V\in\PP_q$ and the function
\[
 U\mapsto P_V^{(\alpha)}F_M(U)
\]
is nondecreasing in the Loewner order for every
$V>0$, $M\geq0$ and every choice of positive
exponents $\nu_2,\ldots,\nu_m$.
In particular, determine the optimal lower bound
on $\alpha$ in the continuous shape range.
\end{question}

In the following two questions, we consider
the real Wishart distribution at discrete shapes
$n=k/2$ and write $D_j=\det S_{jj}$.
\begin{question}\label{ques2}
For prescribed block sizes $p_1,\ldots,p_m$ and
positive exponents $\nu_2,\ldots,\nu_m$, determine
all integers $k\geq\max_j p_j$ for which the
conditional moment
\[
 a\mapsto\E\left[
 \prod_{j=2}^m D_j^{\nu_j}\,\middle|\,S_{11}=a
 \right]
\]
admits a continuous Loewner-nondecreasing version
on the positive definite cone for every $\Sigma>0$.
\end{question}

\begin{question}\label{ques3}
For prescribed block sizes $p_1,\ldots,p_m$ and
positive exponents $\nu_1,\ldots,\nu_m$, determine
all integers $k\geq\max_j p_j$ for which there exists
a permutation $\pi$ of $\{1,\ldots,m\}$, independent
of $\Sigma$, such that
\[
 \E\prod_{j=i}^m D_{\pi(j)}^{\nu_{\pi(j)}}
 \geq
 \E D_{\pi(i)}^{\nu_{\pi(i)}}
 \E\prod_{j=i+1}^m D_{\pi(j)}^{\nu_{\pi(j)}}
\]
for every $\Sigma>0$ and $i=1,\ldots,m-1$.
\end{question}

The properties in Questions~\ref{ques2} and \ref{ques3} hold for every $k\geq p-1$.
When $m=2$, Anderson's inequality extends both
properties to the entire range $k\geq\max_j p_j$.

\section*{Declaration on the use of artificial intelligence}
The authors use AI tools to generate draft proofs of WPI at the initial stage. Subsequent abstraction, refinement, and generalization of the AI's ideas are done by humans. The authors take full responsibility for the use of AI and for the mathematical content, accuracy, and integrity of the manuscript.

\appendix

\section{Transition functions and Markov processes}\label{app:markov-process}

\subsection{Transition functions}\label{sec:markov-process}  Let $E$ be a metric space endowed with its Borel
$\sigma$-algebra $\mathscr{B}$. A function $P(t,x,\varGamma)\ (t\ge 0,x\in E,\varGamma\in\mathscr{B})$ is called a {transition function} if the following conditions are satisfied \cite[Chapter~II, Section~1]{Dynkin1965}
\begin{enumerate}
\item[(M1)] for fixed \(t\) and \(x\), the function $\varGamma\mapsto P(t,x,\varGamma)$ is a measure on \(\mathscr{B}\),
\item[(M2)] for fixed \(t\) and \(\varGamma\), the function
$x\mapsto P(t,x,\varGamma)$ is \(\mathscr{B}\)-measurable,
\item[(M3)] \(P(t,x,E)\le 1\),
\item[(M4)] \(P(0,x,E\setminus\{x\})=0\), and
\item[(M5)] $P(s+t,x,\varGamma)= \int_E P(s,x,\mathrm{d}y)\,P(t,y,\varGamma)\ (s,t\ge 0)$.
\end{enumerate}

From (M3) and (M5) we have
\[
P(s+t,x,E)
\le
\int_E P(s,x,\mathrm{d}y)
=
P(s,x,E) \ (s,t\ge 0).
\]
Consequently, \(P(t,x,E)\) is a nonincreasing function of \(t\), and the limit $P(+0,x,E)$ exists. We say that the transition function \(P(t,x,\Gamma)\) is {normal} if $P(+0,x,E)=1$ for all $x\in E$. We call it {conservative} if $P(t,x,E)=1$ {for all} $t>0$, $x\in E$. By (M3)-(M5), every normal transition function satisfies
$P(0,x,\varGamma)=1_{\varGamma}(x)$.

Let $B_b(E)$ be the Banach space of bounded
$\mathscr B$-measurable real-valued functions on $E$,
equipped with the supremum norm
$\|f\|_\infty=\sup_{x\in E}|f(x)|$.
A family of operators $P_t:B_b(E)\to B_b(E)$
$(t\geq0)$ is called a {contraction
semigroup} if
\begin{enumerate}
\item[(T1)] (Linearity)
$P_t(af+bg)=aP_tf+bP_tg$,
\item[(T2)] (Positivity)
$P_tf\geq0$ whenever $f\geq0$,
\item[(T3)] (Contraction)
$\|P_tf\|_\infty\leq\|f\|_\infty$,
\item[(T4)] (Identity)
$P_0=I$, and
\item[(T5)] (Semigroup property)
$P_sP_t=P_{s+t}$,
\end{enumerate}
for all $s,t\geq0$, $a,b\in\mathbb R$ and
$f,g\in B_b(E)$.

Every normal transition function satisfying (M1)-(M5)
gives rise to a contraction semigroup via
\begin{equation}\label{eq:semigroup}
     P_tf(x)=\int_E f(y)P(t,x,\dd y)\ (f\in B_b(E)).
\end{equation}
The associated operators satisfy the additional condition
\begin{enumerate}
\item[(T6)] (Order continuity)
$P_tf_n\downarrow0$ whenever $f_n\downarrow0$
pointwise.
\end{enumerate}
Conversely, a family of operators satisfying (T1)-(T6)
determines a transition function by
$P(t,x,\varGamma)=P_t1_\varGamma(x)$.
Normality is equivalent to
$\lim_{t\downarrow0}P_t1(x)=1$ for every $x\in E$,
and conservativity is equivalent to $P_t1=1$
for all $t\geq0$. By positivity and linearity, the operators $P_t$
preserve the pointwise order.
They extend to nonnegative measurable functions by
$P_tf=\sup_{n\geq1}P_t(f\wedge n)$.
Under this extension, the semigroup property and
monotone convergence remain valid.

A Markov {kernel} $K(x,\varGamma)\ (x\in E,\varGamma\in\mathscr B)$ satisfies (M1)-(M3) without the time parameter, together with $K(x,E)=1$ for all $x\in E$. Its associated {operator}, also denoted by $K$, acts as $Kf(x)=\int_E f(y)K(x,\dd y)$ for $f\in B_b(E)$, as in~\eqref{eq:semigroup}. Equivalently, a Markov operator on $B_b(E)$ satisfies (T1)-(T3), (T6) and $K1=1$, with $P_t$ replaced by $K$. 

\subsection{Filtered probability spaces}
Let $(\varOmega,\mathscr F,
(\mathscr F_t)_{t\ge0},\mathbb P)$ be a filtered probability space, equipped with a filtration of $\sigma$-subalgebras $\mathscr F_t \ ({t\geq0})$ such that $ \mathscr F_s \subseteq \mathscr F_t \subseteq \mathscr F \ (s\le t)$. We say that the filtration satisfies
the natural conditions 
\cite[Definition~1.3.38]{Bichteler2002} if
\begin{enumerate}
\item[(A1)] ($\sigma$-completeness)
every subset of a $\mathbb P$-null set in $\mathscr F_t$
belongs to $\mathscr F_0$ for each $t<\infty$, and
\item[(A2)] (Right continuity)
$\mathscr F_t=\bigcap_{s>t}\mathscr F_s$
for every $t\geq0$.
\end{enumerate}
The usual conditions are obtained by replacing (A1)
with the requirement that every subset of a
$\mathbb P$-null set in $\mathscr F$ belongs to
$\mathscr F_0$, which is generally stronger than (A1). 

An $E$-valued stochastic process is a family of measurable mappings
$X_t:\varOmega \to E \ ({t\geq0})$. It is adapted if $X_t$ is $\mathscr F_t$-measurable for every $t\geq0$. Let $(\mathbb P_x)_{x\in E}$ be a family of probability measures such that $\mathbb P_x(X_0=x)=1$.
A stopping time is a measurable function $\tau:\Omega\to [0,\infty]$ on the extended real line  such that $\{\omega:\tau(\omega)\le t\} \in \mathscr F_t$ for every $t \ge 0$. An adapted process $X_t  \ ({t\geq0})$ is called a Markov process with
transition semigroup $P_t\ ({t\geq0})$ if
\[
 \E_x[f(X_{s+t})\mid\mathscr F_s]
 =P_tf(X_s)\ (s,t\geq0),
\]
for every $x\in E$ and $f\in B_b(E)$.
It is called a strong Markov process if, for every almost surely finite stopping time $\tau$,
\[
 \E_x[f(X_{\tau+t})\mid\mathscr F_\tau]
 =P_tf(X_\tau) \  (t\geq0),
\]
$\mathbb P_x$-almost surely for every $x\in E$
and $f\in B_b(E)$.

\section{Wishart diffusions on symmetric cones}\label{sec:wishart-semigroup}
Let $\mathcal C$ be a closed convex cone in a
finite-dimensional real vector space. Suppose that
$P_A\ (A\in\mathcal C)$ is a family of positive linear
contraction operators on $B_b(\mathcal C)$ satisfying
(T6), $P_A1=1$, and the following two conditions
\begin{enumerate}
\item[(T4$'$)] (Identity) $P_0=I$, and
\item[(T5$'$)] (Semigroup property) $P_AP_B=P_{A+B}$ ($A,B\in\mathcal C$).
\end{enumerate}
For every nonzero $H\in\mathcal C$, the family
$P_{tH}\ (t\geq0)$ is a conservative transition
semigroup. The set $\varGamma = \{x H:x \geq0\}$ is called the ray generated
by $H$. It is invariant under $P_{tH}\ (t\geq0)$ if $P_{tH} 1_{\varGamma}(xH)=1 \ (t,x\geq0).$

Consider an important class of affine diffusions. Fix $\alpha>0$ with $2\alpha\geq\beta(q-1)$. For $T>0$ and
$U,V\in\PP_q$, put
\[
 e_T(U)=\exp\{-\tr(TU)\} \quad \text{and} \quad \psi_V(T)=(T^{-1}+2V)^{-1}.
\]

\begin{lemma}
\label{prop:wishart-semigroup}
For every $V\in\PP_q$, there exists a unique probability
kernel $P_V^{(\alpha)}(U,\varGamma)$ on $\PP_q$ characterized
by the affine transform
\begin{equation}\label{eq:cond-transform}
 P_V^{(\alpha)}e_T(U)
 =\left(\frac{\det\psi_V(T)}{\det T}\right)^\alpha
 e_{\psi_V(T)}(U).
\end{equation}
The associated operators satisfy (T4$'$) and (T5$'$).
Moreover, for every $V\in\PP_q$, the transition semigroup
$P_{tV}^{(\alpha)}\ (t\geq0)$ admits a continuous
Markov realization.
\end{lemma}

\begin{proof}[Proof of Lemma \ref{prop:wishart-semigroup}]
The existence of the kernels and of continuous
Markov realizations of the associated semigroups
follow from \cite[Section~5.1, Propositions~5.5 and~5.7]{CKMT2016}.
Uniqueness follows from the uniqueness of Laplace
transforms.

Suppressing the superscript $\alpha$, we obtain
$P_0=I$ from~\eqref{eq:cond-transform}.
Moreover, the identity
$\psi_V\circ\psi_H=\psi_{V+H}$ and the cancellation
of the determinant ratios give $P_VP_He_T=P_{V+H}e_T.$
Uniqueness of Laplace transforms yields
\begin{equation}\label{eq:cond-addition}
 P_VP_H=P_{V+H}=P_HP_V\ (V,H\in\PP_q).
\end{equation}
Thus \textnormal{(T4$'$)} and \textnormal{(T5$'$)} are satisfied. In particular, $P_{V}(0,\varGamma)$
is the central Wishart distribution
$\W_q^{(\beta)}(\alpha,\beta V)$ evaluated at $\varGamma$.
\end{proof}

The kernel $P_V^{(\alpha)}$ determined by Lemma~\ref{prop:wishart-semigroup} is called a Wishart kernel with shape parameter $\alpha$ and covariance parameter $\beta V$. 

A function $f:\PP_q\to\R$ is said to be of {polynomial growth} if there exist constants
$C>0$ and $h\geq0$ such that
\[
 |f(U)|\leq C(1+\tr U)^h\ (U\in\PP_q).
\]
Since the trace is equivalent to any matrix norm 
on $\PP_q$, this definition is independent of the
choice of norm.

\begin{lemma}\label{lem:wishart-continuity}
For every positive integer $h$,
\begin{equation}\label{eq:cond-moments}
 \int_{\PP_q}(1+\tr X)^h P_V(U,\dd X)
 \leq C_h(1+\tr U)^h,
\end{equation}
where $C_h$ can be chosen uniformly for $V$ in
compact subsets of $\PP_q$.
Consequently, if $f$ is continuous and of polynomial
growth, then $P_Vf(U)$ is finite and jointly continuous
in $(V,U)$ on $\PP_q\times\PP_q$.
\end{lemma}

\begin{proof}[Proof of Lemma \ref{lem:wishart-continuity}]
Taking $T=sI$ in~\eqref{eq:cond-transform}, we obtain
\[
 P_Ve_{sI}(U)
 =\det(I+2sV)^{-\alpha}
 \exp\{-s\tr((I+2sV)^{-1}U)\}.
\]
For $k=0,\ldots,h$, the $k$th moment of $\tr X$ is
\[
 \int_{\PP_q}(\tr X)^k P_V(U,\dd X)
 =(-1)^k\left.\frac{\dd^k}{\dd s^k}P_Ve_{sI}(U)\right|_{s=0+}.
\]
It is a polynomial in $U$ of degree at most $k$,
with coefficients continuous in $V$.
This gives~\eqref{eq:cond-moments}.

The continuity of the Laplace transform in $(V,U)$
implies weak continuity of the kernels.
The moment bound~\eqref{eq:cond-moments}, applied
with a sufficiently large $h$, ensures uniform
integrability along convergent sequences of $(V,U)$.
The asserted continuity follows.
\end{proof}

Fix a direction $H\in\PP_q$. Let $U_t \ (t \ge 0)$ be the continuous Wishart diffusion
with semigroup $P_{tH}^{(\alpha)} \ ({t\geq0})$, and write $\mathbb P_a$
for its law with $U_0=a$. 
\begin{lemma}[Local characteristics]\label{lem:wishart-characteristics}
Let $U_t\ (t\geq0)$ be a continuous Wishart diffusion
with transition semigroup $P_{tH}\ (t\geq0)$.
For a Hermitian matrix $T$, the process
\[
 N_t^T=\tr(TU_t)-\tr(TU_0)-2\alpha t\tr(TH)
\]
is a continuous martingale with quadratic variation
\begin{equation}\label{eq:cond-bracket}
 \dd\langle N^T\rangle_t=4\tr(U_tTHT)\,\dd t.
\end{equation}
\end{lemma}

\begin{proof}[Proof of Lemma \ref{lem:wishart-characteristics}]
Differentiating~\eqref{eq:cond-transform} in the
Laplace variable gives
\[
 P_{tH}\bigl[X\mapsto\tr(TX)\bigr](U)=\tr(TU)+2\alpha t\tr(TH).
\]
Moreover, the variance of $\tr(TU_t)$ under $\mathbb P_U$
is
\[
 4t\tr(UTHT)+4\alpha t^2\tr(THTH).
\]
Both identities extend to all Hermitian $T$ because both sides of each identity are polynomials in the real coordinates of $T$.
By the Markov property, $N_t^T$ is a martingale.
The conditional second moment of its increments,
together with the first moment formula, shows that
\[
 (N_t^T)^2-4\int_0^t\tr(U_sTHT)\,\dd s
\]
is also a martingale. This proves~\eqref{eq:cond-bracket}.
\end{proof}


\bibliographystyle{amsplain}
\bibliography{bibtex}
\end{document}